\documentclass[11pt,oneside]{amsart}

\usepackage{blindtext}

\usepackage{amsbsy,amsmath,amssymb}

\def\E{\mathbb{E}}

\def\R{\mathbb{R}}

\def\N{\mathbb{N}}

\def\P{\mathbb{P}}

\def\cA{\mathcal{A}}
\def\cB{\mathcal{B}}

\def\cD{\mathcal{D}}

\def\cE{\mathcal{E}}

\def\cH{\mathcal{H}}

\def\cL{\mathcal{L}}

\def\cN{\mathcal{N}}

\def\ee{\boldsymbol{e}}

\def\gg{\boldsymbol{g}}

\def\mm{{\boldsymbol{m}}}

\def\sb{\boldsymbol{s}}

\def\uu{\boldsymbol{u}}
\def\vv{\boldsymbol{v}}
\def\ww{\boldsymbol{w}}

\def\xx{\boldsymbol{x}}
\def\yy{\boldsymbol{y}}
\def\zz{\boldsymbol{z}}

\def\BB{\boldsymbol{B}}

\def\GG{\boldsymbol{G}}
\def\HH{\boldsymbol{H}}
\def\II{\boldsymbol{I}}

\def\KK{\boldsymbol{K}}

\def\QQ{\boldsymbol{Q}}

\def\YY{\boldsymbol{Y}}

\def\00{\boldsymbol{0}}
\def\11{\boldsymbol{1}}
\def\22{\boldsymbol{2}}
\def\33{\boldsymbol{3}}
\def\44{\boldsymbol{4}}
\def\55{\boldsymbol{5}}
\def\66{\boldsymbol{6}}
\def\77{\boldsymbol{7}}
\def\88{\boldsymbol{8}}
\def\99{\boldsymbol{9}}

\def\aalpha{{\boldsymbol{\alpha}}}
\def\bbeta {{\boldsymbol{\beta}}}

\def\eepsilon{{\boldsymbol{\epsilon}}}

\def\pphi{{\boldsymbol{\phi}}}

\def\ssigma{{\boldsymbol{\sigma}}}

\def\oomega{{\boldsymbol{\omega}}}

\newcommand{\bs}[1]{\boldsymbol{#1}}
\newcommand{\cEE}{\bs{\mathcal{E}}}

\def\PPhi{{\boldsymbol{\Phi}}}

\DeclareMathOperator*{\spn}{span}
\DeclareMathOperator{\supp}{supp}
\def\transpose{{\rm T}}
\def\curl{\qopname\relax o{curl}}
\def\div{\qopname\relax o{div}}

\usepackage{xcolor}
\newcommand{\tlg}[1]{{\color{black}{#1}}}
\newcommand{\is}[1]{{\color{black}{#1}}}
\newcommand{\hw}[1]{{\color{black}{#1}}}

\newtheorem{theorem}{Theorem}[section]
\newtheorem{corollary}[theorem]{Corollary}

\newtheorem{lemma}[theorem]{Lemma}
\newtheorem{proposition}[theorem]{Proposition}
\newtheorem{definition}[theorem]{Definition}

\newtheorem{remark}[theorem]{Remark}

\newtheorem{assumption}[theorem]{Assumption}

\usepackage{lastpage}

\begin{document}

\title[Vector-Valued Gaussian Processes]{Vector-Valued Gaussian Processes for Approximating Divergence-
  or Rotation-free Vector Fields}

\author{Quoc Thong Le Gia}
\address{Quoc Thong Le Gia, School of Mathematics and Statistics,
  University of New South Wales, Sydney, NSW 2052, Australia}
\email{qlegia@unsw.edu.au}

\author{Ian Hugh Sloan}
\address{Ian Hugh Sloan, School of Mathematics and Statistics,
  University of New South Wales,
  Sydney, NSW 2052, Australia}
\email{i.sloan@unsw.edu.au}

\author{Holger Wendland}
\address{Chair of Numerical and Applied Mathematics, Department of
  Mathematics, University of Bayreuth, 95440 Bayreuth, Germany}
\email{holger.wendland@uni-bayreuth.de}

\begin{abstract}%
\noindent In this paper, we discuss vector-valued Gaussian processes for the
approximation of di\-ver\-gence- or rotation-free functions. We establish
the theory for such Gaussian processes, then link the theory to
multivariate approximation theory, and finally give error estimates for the
predictive mean in various situations.
\end{abstract}

\maketitle

\section{Introduction}

Gaussian processes are established tools for approximating unknown
functions under the influence of uncertainties
\cite{Rasmussen-Williams-06-1}. The approximation 
properties of Gaussian processes of scalar-valued functions, the {\em
  so-called regression problem}, have 
extensively been studied, see for example
\cite{Nieman-etal-22-1,Pati-etal-15-1,Teckentrup-20-1,
  Vaart-etal-11-1,Wynne-etal-21-1,Zhou-etal-23-1}.
However, many applications, particularly from
physics and engineering, require the reconstruction of vector-valued
or, {more} generally, Hilbert-space valued functions. For example, when
modeling fluid flow problems with incompressible flows, the velocity
field at a fixed point in time is given by a differentiable function $\vv:\cD\to\R^d$
with  $\cD\subseteq\R^d$ and usually $d=2,3$ satisfying $\nabla\cdot\vv =
\sum_j \partial_j v_j =0$ if $\vv=(v_1,\ldots,v_d)$.

Another application is the modeling of magnetic fields using
Maxwell's equation{s}. Here two fields $\HH:\cD\to\R^3$ and
$\BB:\cD\to{\R^3}$ are usually used, which have to satisfy, among other things,
$\nabla\cdot\BB=0$ and $\nabla\times \HH=\00$, respectively.

In this paper, we do not {consider the underlying}
partial differential equations, but rather assume that we have some
{\em training data} in the form of {\em data sites}
$X=\{\xx_1,\ldots,\xx_N\}\subseteq\cD$ and {vector} {\em data values}
$\vv(\xx_j)$, $1\le j\le N$, with $\vv: \cD \to \R^3$.

It is the goal of this paper to define precise generalizations of
scalar-valued Gaussian processes for vector-valued or even
Hilbert-space-valued functions. We are particularly interested in
modeling vector-valued functions which are either divergence or
rotation free.

While there are already some papers that address at least the
modeling question, see for example \cite{Parra-Tobar-17-1, Bilionis-etal-13-1} for 
general vector-valued functions and \cite{Wahlstrom-etal-13-1} and {\cite{Owhadi-23-1}}  for
divergence- and rotation-free functions, there is, so far, no error
analysis of such processes. Only very first preliminary steps in this
direction can be found in \cite{Li-etal-??-1}, though they differ
significantly from ours.

The new contributions of this paper are the mathematically precise
introduction and discussion of divergence- and rotation-free Gaussian
processes {together with} the analysis of their connection to matrix-valued positive
definite kernels and vector-valued reproducing kernel Hilbert spaces,
and {finally} error analysis in Sobolev spaces. For the latter, we can rely
on results from approximation theory, particularly
\cite{Fuselier-08-1}, but {require} new results on the approximation
operator. 

\section{Gaussian Processes and Approximation}
  
We shortly review the necessary material on (scalar-valued) Gaussian
processes and their connection to multivariate, kernel-based
approximation theory, laying out the path we will also follow for
vector-valued functions.

\subsection{Scalar-Valued Gaussian Processes}
According to
\cite{Rasmussen-Williams-06-1}: ``A Gaussian process is a collection of
random variables, any finite number of which have a joint Gaussian
distribution''.
To be more precise, assume initially that we have a probability space
$(\Omega,\cA,\P)$ consisting of the abstract domain
\is{$\Omega$}, a $\sigma$-\is{algebra} $\cA$ on $\Omega$ and a
probability measure $\P$ on $(\Omega,\cA)$.

Given an {\em index domain} $\cD\subseteq\R^d$, a Gaussian process is
a family of random variables $g:\Omega\times\cD\to\R$, such that for 
any finite set of pairwise distinct points 
$\{\xx_1,\ldots,\xx_N\}\in\cD$, the random variable $\gg=(g(\cdot,\xx_1),
\ldots, g(\cdot,\xx_N))^\transpose$  is jointly \is{normally} distributed,
i.e. it possesses the probability density function
$\mathfrak{g}:\R^N\to\R$ of the form
\[
\mathfrak{g}(\zz) =  (2\pi)^{-N/2} \det(\Sigma)^{-1/2}
\exp\left(-\frac{1}{2}(\zz-\mm)^\transpose\Sigma^{-1}
(\zz-\mm)\right), \qquad \zz\in \R^N,
\]
where $\mm=(m(\xx_1),\ldots, m(\xx_N))^\transpose$ is the mean,
and the positive definite and symmetric matrix $\Sigma =
(\E[(g(\cdot,\xx_i)- m(\xx_i))(g(\cdot,\xx_j)-m(\xx_j))])_{1\le 
  i,j\le N} \in \R^{N\times N}$ is the covariance matrix. Obviously, this
means that the Gaussian process is uniquely determined by a function
$m:\cD\to\R$ giving the mean  and  a function $K:\cD\times\cD\to\R$
describing the covariance via $K(\xx_i,\xx_j)=\Sigma_{ij}$.
As usual in this context, we will
omit the dependence on the unknown variable from $\Omega$ and simply
write $g(\xx)$ instead of $g(\cdot,\xx)$.  \is{Throughout we shall} use the notation
\[
g\sim {\rm GP}(m, K).
\]
The Gaussian process is called {\em centered} if it has mean zero.

Gaussian processes are often used to model the approximation of an
unknown function $f:\cD\to\R$. If nothing at all is known about this
function, then we may try to approximate $f$ by \is{an initial} Gaussian
process $f_0:\Omega\times \cD\to\R$ with mean function $m(\xx)=0$, $\xx\in\cD$,
though other mean functions are possible, and
a preselected kernel function $K:\cD\times\cD\to\R$. Often, this
kernel depends on several parameters, but we
will suppress this for simplicity.  Obviously, $f_0$ \is{may} not be a particularly good
approximation.  \is{This will} depend significantly on the choice of the
kernel.
\\
\is{Importantly}, if we have $N$ observations of the unknown function $f$ at
data sites $X=\{\xx_1,\ldots,\xx_N\}\subseteq\cD$ \is{(also called {\em
  design or training points)}, with data values either in the form $y_j=f(\xx_j)$ $1\le
j\le N$ without noise, or in 
the form $y_j=f(\xx_j)+\epsilon_j$, $1\le j\le N$ with noise,} than we
can condition the initial Gaussian process using these data to obtain
an {\em a \is{posteriori} Gaussian process}. In the case of no noise,
$f(\xx_j)=y_j$, this process is
given by
$f_{N,0} \sim {\rm GP} (m_N,K_N)$ with new mean and kernel functions
\begin{eqnarray}
m_N(\xx) &:= &m(\xx)  + K(\xx,X) K(X,X)^{-1}\left(\yy- m(X)\right),\label{eq:Gauss1}\\
K_N(\xx,\xx')&:=& K(\xx, \xx') - K(\xx,X) K(X,X)^{-1}K(\xx',X)^\transpose,\label{eq:Gauss2}
\end{eqnarray}
where $\yy=(y_1,\ldots, y_N)^\transpose \in\R^N$ are the observations
at 
the design points
$X$, $m(X):= (m(\xx_1), \ldots, m(\xx_N))^\transpose \in
\R^N$ is the vector of means and
 $K(\xx, X) := (K(\xx,\xx_1), \ldots, K(\xx, \xx_N)) \in \R^{1\times N}$.  \is{Finally,} $K(X,X) \in \R^{N\times N}$
is the matrix with elements $K(\xx_i, \xx_j)$.

The functions $m_N$ and $K_N$ are also called {\em predictive mean}
and {\em predictive covariance}, \is{respectively}. 

In the context of noisy observations $y_j=f(\xx_j)+\epsilon_j$, $1\le
j\le N$, the assumption that the noise $\epsilon_j$, $1\le j\le N$,
follows an independent, identical Gaussian 
distribution with mean zero and variance $\sigma^2>0$ leads to  the \is{
posterior} Gaussian process  $f_{N,\sigma}\sim {\rm
  GP}(m_{N,\sigma}, K_{N,\sigma})$ with 
\begin{eqnarray}
m_{N,\sigma}(\xx) &:=& m(\xx)  + K(\xx,X)(K(X,X)+\sigma^2 I)^{-1}\left(\yy- m(X)\right),\label{eq:Gauss3}\\
K_{N,\sigma}(\xx,\xx')&:= &K(\xx, \xx') - K(\xx,X) (K(X,X)+\sigma^2
I)^{-1}K(\xx',X)^\transpose \label{eq:Gauss4} .
\end{eqnarray}
Here, $I\in\R^{N\times N}$ denotes the identity matrix. Otherwise, the
notation is as in (\ref{eq:Gauss1}) and (\ref{eq:Gauss2}).
From the above equations, we see that $f_{N,\sigma}$ reduces to
$f_{N,0}$ if  \is{$\sigma=0$, justifying} our notation.

\subsection{Connection to Multivariate Approximation Theory}
There is a strong connection to multivariate interpolation theory,
particularly approximation by kernels or {\em radial basis functions.}
To specify this connection,  we will require the following terminology.

\begin{definition}
\begin{enumerate}
\item A continuous, symmetric function  $K:\cD\times\cD\to\R$ is called
{\em positive definite} if for all $N\in\N$ and all sets
points, the kernel matrix $K(X,X)\in\R^{N\times N}$ is positive
definite (that is, if $\hw{\aalpha^\transpose K(X,X) \aalpha} > 0$ for all
$\hw{\aalpha} \in \R^N\setminus \tlg{\{\mathbf{0}\}}$). Such a
  function is also called a {\em positive definite kernel}. 

\item A Hilbert space $\cH=\{f:\cD\to\R\}\subseteq C(\cD)$ of
  functions is a {\em reproducing kernel Hilbert space (RKHS)} if there
    is a function $K:\cD\times \cD\to\R$ satisfying $K(\cdot,\xx)\in
    \cH$ for all $\xx\in \cD$ and
    $
    f(\xx) = \langle f,K(\cdot,\xx)\rangle_\cH
    $
    for all $f\in \cH$ and all $\xx\in\cD$. The function $K$ is called
    the {\em reproducing kernel} of $\cH$.
  \end{enumerate}
  \end{definition}
It is well known, see for example \cite{Aronszajn-50-1,
  Wendland-05-1}, that the reproducing kernel of a RKHS is always
positive semi-definite (i.e. $K(X,X)$ is symmetric and positive
semi-definite for all $X$) and that the kernel is uniquely determined
by the Hilbert space. Moreover, for every positive definite kernel
there is a unique RKHS having this kernel as its reproducing kernel.

For us, the following connections between positive definite kernels,
RKHS and Gaussian processes are important. We start with the case of
no noise.

\begin{theorem} Let $K:\cD\times\cD\to\R$ be a positive definite
  kernel and let $X=\{\xx_1,\ldots,\xx_N\}\subseteq\cD$ consist of
  pairwise distinct points. Define the  \is{$N$-}dimensional
  approximation space  $V_X=\spn\{K(\cdot,\xx_j) : 1\le j\le N\}$ and
  let $\cH$ be the RKHS of the kernel $K$.
  \begin{enumerate}
  \item For every $\yy\in \R^N$ there is a unique function $s_{0,\yy}\in
    V_X$ interpolating the data,
      i.e. satisfying $s_{0,\tlg{\yy}}(\xx_i) = y_i$, $1\le i\le N$. This defines
      an interpolation operator $I_X:\R^N\to V_X$, $I_X\yy:=s_{0,\yy}$
      and also an operator $I_X:C(\cD)\to V_X$ by setting
      $I_Xf=I_X(f|X)$.
    \item The function $s_{0,\yy}\in V_X$ is the solution of
      \[
      \min\left\{\|s\|_\cH : s\in \cH \mbox{ with } s(\xx_i) = y_i,
      1\le i\le N\right\}.
      \]
    \item If $m_N$ is the predictive mean from (\ref{eq:Gauss1}) with
      data given by $y_i=f(\xx_i)$, $1\le i\le N$,  then  $m_N=m+I_X
      f-I_X m$.
      \item If $K_N$ is the predictive covariance kernel from
        (\ref{eq:Gauss2}) then
        \[
        K_N(\xx,\xx') = K(\xx,\xx') -I_X K(\cdot,\xx') (\xx) =
        K(\xx,\xx') - I_X K(\cdot,\xx) (\xx'),
        \]
        i.e. the predictive covariance kernel is the error between the
        given kernel and its interpolant with respect to one of its arguments.

        For the predictive variance $\sigma_N^2(\xx):=K_N(\xx,\xx)$
        this means
        \begin{align*}
          \sigma_N^2(\xx) & =  \sup_{\|f\|_{\cH} =1} |f(\xx) -I_Xf
          (\xx)|^2 \\
          & =  K(\xx,\xx)-2K(\xx,X) \uu(\xx) +
          \uu(\xx)^\transpose K(X,X)\uu(\xx),
        \end{align*}
        where the vector
        $\uu(\xx)=(u_1(\xx),\ldots, u_N(\xx))^\transpose\in\R^N$
        contains the cardinal or Lagrange functions 
        $u_j\in V_X$ evaluated at $\xx$. They satisfy
        $u_j(\xx_i)=\delta_{ij}$. 
  \end{enumerate}
\end{theorem}
\begin{proof}
  The first two statements are well known in the approximation
  literature, see for example \cite{Wendland-05-1}.\is{ An explicit formula for the interpolant follows from writing the interpolant in the form}
  \[
  s_{0,\yy} = \sum_{j=1}^N \alpha_j K(\cdot,\xx_j),
  \]
  \is{and then using the interpolation conditions to give $\aalpha = K(X,X)^{-1} \yy$,
  and hence 
  \begin{equation}\label{interp}
  (I_f \,\yy)(\xx) = s_{0,\yy}(\xx) = K(\xx,X) \aalpha =
  K(\xx,X) K(X,X)^{-1}\yy.
  \end{equation}} 
\is{The \tlg{third} statement follows easily.}
These formulas also immediately imply the first part of the fourth point, concerning the representation of the predictive covariance.

The first representation of the predictive variance is also known in the
approximation community for quite some time, see 
\cite{Scheurer-etal-13-1} and \cite{Stuart-Teckentrup-17-1}, for example.
The remaining representation follows from the fact that the cardinal
functions are defined by the equation $K(X,X)\uu(\xx)=K(\xx,X)^\transpose$.
\end{proof}

In the approximation theory literature the expression
$\sigma_N^2(\xx)$ is usually called the {\em Power function} of the
approximation process and is denoted by $P_{K,X}^2(\xx)$.

The connection between Gaussian processes and RKHS and positive
definite kernels in the situation of noisy data is summarized in the
next theorem. It is a specific version of the {\em representer
  theorem} from machine learning, see
\cite{Schoelkopf-Smola-02-1,Steinwart-Christmann-08-1}. 

\begin{theorem}\label{thm:representer}
  Let $K:\cD\times\cD\to\R$ be a positive definite
  kernel and let $X=\{\xx_1,\ldots,\xx_N\}\subseteq\cD$ consist of
  pairwise distinct points. Define the \is{$N$-} dimensional
  approximation space  $V_X=\spn\{K(\cdot,\xx_j) : 1\le j\le
  N\}$ and let $\cH$ be the RKHS of the kernel $K$.
  \begin{enumerate}
  \item Given $\lambda>0$ and $\yy\in\R^N$, there is a unique function $s_{\lambda,\yy}\in V_X$
    solving
    \[
    \min\left\{\sum_{j=1}^N |s(\xx_j)-y_j|^2 + \lambda\|s\|_\cH^2 :
    s\in\cH\right\}.
    \]
  This defines an approximation operator $Q_{X,\lambda}:\R^N\to V_X$
  $Q_{X,\lambda}\yy:=s_{\lambda,\yy}$ and also an operator $Q_{X,\lambda}:C(\cD)\to V_X$  
  by setting $Q_{X,\lambda}f := Q_{X,\lambda} (f|X)$.
  \item
  If $m_{N,\sigma}$ is the predictive mean from (\ref{eq:Gauss3}) with
  data $\yy\in\R^N$ and   prior mean $m$ and covariance kernel $K$, then $m_{N,\sigma} = m-
  Q_{X,\sigma^2}m + Q_{X,\sigma^2}\yy$.
  \end{enumerate}
\end{theorem}

Again, we see the connection between the no-noise and the noisy case,
as we have $Q_{X,0} = I_X$. Moreover, we can use the approximation
operator also for no-noise data\is{, and the interpolation operator for noisy data with either}
known or unknown distribution. 

The connection between the mean of the predictive or \is{posterior}
Gaussian process and kernel based interpolation and approximation,
\is{leads to error estimates} typically of the form
\[
\|f(\xx)- m_N(\xx)\|_{L_\infty(\cD)} \le C h_{X,\cD}^{\tau-d/2}
\|f\|_{\cH},
\]
where $h_{X,\cD}$ is the so-called fill distance of $X$ in $\cD$,
defined as
\begin{equation}\label{filldistance}
h_{X,\cD} = \sup_{\xx\in\cD} \min_{1\le j\le N} \|\xx-\xx_j\|_2.
\end{equation}
\is{Here the data-giving} function $f$ is assumed to be in the Sobolev space
$\cH=H^\tau(\cD)$, $\cD$ is supposed to have a 
Lipschitz boundary and $K$ is supposed to be the reproducing kernel of
$H^\tau(\cD)$. More general estimates that involve the error and
relax the condition on the smoothness of $f$ can be found in
    \cite{Wynne-etal-21-1}. \is{Moreover}, as we also have identified the
variance as the power function, we can use bounds on the power
function and classical estimates like the Chebyshev inequality to
derive estimates on the whole process, not only the mean. A prototype
of such an estimate is the following
\is{bound on the departure of the interpolation estimate from the predictive mean.}
\begin{theorem} Assume that the covariance kernel $K$ is the
  reproducing kernel of a Sobolev space $H^\tau(\cD)$ with $\tau>d/2$. Then, there
  are constants $C>0$ and $h_0>0$ such that for all 
  $X=\{\xx_1,\ldots,\xx_N\}\subseteq\cD$ with $h_{X,\cD}\le h_0$, all
  $f\in H^\tau(\cD)$ and
  every $\epsilon \in (0,1)$, the bound
  \[
  \P\left(|f_{N,0}(\xx)-m_N(\xx)|\ge\epsilon\right) \le C\frac{h^{2\tau-d}_{X,\cD}}{\epsilon^2}
  \]
  holds.
\end{theorem}
\begin{proof}
  This follows from Chebyshev's inequality
  \[
  \P\left(|f_{N,0}(\xx)-m_N(\xx)|\ge\epsilon\right) \le
  \frac{\sigma_N(\xx)^2}{\epsilon^2}
  \]
  and
  \[
  \sigma_N(\xx)^2 = K_N(\xx,\xx)  = P_{X,K}^2(\xx) \le C
  h_{X,\cD}^{2\tau-d},
  \]
  where the latter is a standard estimate from radial basis function approximation.
\end{proof}
\subsection{Connection to Functional Analysis}

There is yet another thing we need to recall\is{, namely} the connection of
positive definite kernels to the spectral theory of positive,
self-adjoint and compact integral operators. If $K:\cD\times\cD\to\R$
is a positive definite kernel then $A:L_2(\cD)\to L_2(\cD)$ defined by
\[
A f = \int_{\cD} f(\yy) K(\cdot,\yy)d\yy, \qquad f\in L_2(\cD),
\]
has a spectral decomposition, i.e. there is \tlg{an} orthonormal basis of
$L_2(\cD)$ consisting of eigenfunctions of $A$\is{, for which the corresponding  eigenvalues are all}
non-negative. As we only require the positive eigenvalues let
$\{\phi_k\}\subseteq L_2(\cD)$ be an orthonormal set of eigenfunctions \is{corresponding}
to the eigenvalues  
$\lambda_k> 0$, $k\in\N$. Thus, we have $\langle
\phi_k,\phi_j\rangle_{L_2(\cD)} = \delta_{jk}$, 
$A\phi_k=\lambda_k\phi_k$ and
\[
Af = \sum_{j=1}^\infty \lambda_j \langle f,\phi_j\rangle_{L_2(\cD)} \phi_j, \qquad
f\in L_2(\cD).
\]
This can be used to express not only the kernel in the orthonormal
basis but also any Gaussian process. 

\begin{theorem}[Mercer]\label{thm:mercer}
Let $\cD$ be a compact subset of $\R^d$ and let $K:\cD \times \cD \to
\R$ be a positive definite kernel. Let  $\{\phi_k\}_{k\in\N}$ be
an orthonormal subset of $L_2(\cD)$ consisting of eigenfunctions  of
the integral operator with kernel $K$ with corresponding positive
eigenvalues $\{\lambda_k\}_{k\in\N}$. 
Then
\[
K(\xx,\xx') = \sum_{k=1}^\infty 
\lambda_k \phi_k(\xx)\phi_k(\xx'),
\quad \xx, \xx' \in \cD,
\]
where the series converges absolutely and uniformly on $\cD\times \cD$.
\end{theorem}
\is{We can then} express stochastic processes \is{in terms of the \tlg{eigenfunctions} $\phi_k$, as is} done here for Gaussian processes.

\begin{theorem}[Karhunen-Lo\`{e}ve] 
Let $\cD\subseteq\R^d$ be compact. Let $g:\Omega\times\cD\to\R$ be an
integrable, centered Gaussian process with continuous covariance
$K$. Let  $\{\phi_k\}_{k\in\N}$ and  $\{\lambda_k\}_{k\in\N}$ be as in
Theorem~\ref{thm:mercer}. Then, for all $\xx\in\cD$, $\omega\in\Omega$,
\begin{equation}\label{eq:KL}
g(\omega, \xx) = \sum_{k=1}^\infty 
\sqrt{\lambda_k} \xi_k(\omega) \phi_k(\xx),
\end{equation}
where the convergence is in the $L_2(\Omega)$ sense and is uniform on $\cD$.
The random variables $\xi_k$ are independent standard normal random variables. They are given \is{in terms of $g(\omega, \xx)$} by
\[
 \xi_k(\omega) = \frac{1}{\sqrt{\lambda_k}} \int_{\cD} g(\omega, \xx) \phi_k(\xx)
 d\xx.
 \]
\end{theorem}

\section{Vector-Valued Function Approximation}
In this section, we want to generalize the concept of Gaussian
processes to Hilbert-space valued functions. We are particularly
interested in vector-valued Gaussian processes which model physical
vector fields from fluid dynamics and electromagnetics. Thus, we will
introduce \is{divergence-free} or rotation-free Gaussian processes.

After that we will discuss generalizations of positive definite
kernels and reproducing kernel Hilbert spaces to the situation of
Hilbert-space valued functions and derive the connection to Gaussian
processes. 

\subsection{Generalized Gaussian Processes}

The generalization of Gaussian processes to Hilbert-space valued
functions is, more or less, straightforward, \tlg{see also
  \cite[Ch. 1, Sec. 7, p. 35]{Adler-90-1} for separable
  Banach-space-valued Gaussian processes}. 

\is{We} will use the following notation. Let $W$ be a
Hilbert space and let $\cL(W)$ be the space of bounded and linear
operators $A:W\to W$. On $W$ we can, as usual, define the Borel
$\sigma$-algebra $\cB(W)$ which is the smallest $\sigma$-algebra on
$W$ which contains the open sets of $W$. Measurabilty of functions
mapping into $W$ will always be with respect to $(W,\cB(W))$. Next,
for $N\in\N$ we let $W^N=W\times \cdots \times W$ be the $N$-fold Cartesian
product of $W$, which is equipped with the Borel $\sigma$-algebra
$\cB(W^N)=\cB(W)\times \cdots\times \cB(W)$.

\begin{definition}Let  $(\Omega,\cA,\P)$ be a probability space
 and let $W$ be a separable Hilbert space with 
orthonormal basis $\{\ee_j\}_{j\in\N}$. A measurable function
  $\gg:\Omega\to W$ has a {\em Gaussian or normal distribution} if the
  family $(\langle \gg, \ee_k\rangle_W)_{k\in\N}$ is \is{Gaussian in the sense that finitely many of the Fourier coefficients have a joint normal distribution.}
For a given $\cD\subseteq\R^d$, a function $\gg:\Omega\times\cD\to W$
describes a {\em generalized or $W$-valued} Gaussian process if for any finite set of
pairwise distinct points $\{\xx_1,\ldots,\xx_N\}\subseteq\cD$, the
random variable
$\GG:=(\gg(\cdot,\xx_1),\ldots,\gg(\cdot,\xx_N)):\Omega\to W^N$ 
is jointly \is{normally} distributed, i.e. if for any finite index set
$J\subseteq\N$, the family $(\langle \gg(\cdot,\xx_i),
\ee_k\rangle_W)_{1\le i\le N, k\in J}$ has a joint normal distribution.
\end{definition}

In the case $W=\R$ the above definition reduces the above concept of a
normal distribution and a generalized Gaussian process to the
classical concept.

The distribution of a $W$-valued \is{Gaussian random variable} $\gg:\Omega\to W$ is
uniquely determined by its mean $\mm(\gg)$ and covariance matrix
$\Sigma(\gg)$, which take the form
\begin{eqnarray*}
m(\gg)_k &:=& \E[\langle \gg, \ee_k\rangle_W], \qquad k\in\N,\\
\Sigma(\gg)_{k\ell} &:=& \E[(\langle \gg,\ee_k\rangle_W-m(\gg)_k)
  (\langle \gg,\ee_\ell\rangle_W- m(\gg)_\ell)],\ \qquad k,\ell\in\N.
\end{eqnarray*}
The same is true for a $W$-valued Gaussian process \is{$g$}, where mean
$\mm:\cD\to W$ and covariance $\KK:\cD\times\cD\to \cL(W)$ are now determined by
\begin{eqnarray*}
\mm(\xx) &:=& \mm(\gg(\cdot,\xx)), \\
\KK(\xx,\xx')_{k\ell} &:= &\E[(\langle
  \gg(\cdot,\xx),\ee_k\rangle_W-m(\xx)_k)(\langle \gg(\cdot,\xx'), 
  \ee_\ell \rangle_W -m(\xx')_\ell)].
\end{eqnarray*}
We can interpret the kernel indeed as a map $K:\cD\times\cD\to \cL(W)$.
Any element $\vv \in W$ can uniquely be written as $\vv= \sum\langle
\vv,\ee_k\rangle_W \ee_k$ with Fourier coefficients $(\langle
\vv,\ee_k\rangle_W)_{k\in\N}$. Hence, for such an element we can
define  $\tlg{\KK}(\xx,\xx')\vv$ by
\[
\KK(\xx,\xx')\vv : = \sum_{k\in\N}\left( \sum_{\ell\in\N} \tlg{\KK}(\xx,\xx')_{k\ell} \langle
  \vv,\ee_\ell\rangle_W \right)\ee_k.
  \]

We are particularly interested in the case of vector-valued functions
$\gg:\cD\to\R^n$, i.e. where $W=\R^n$, and here our main interest will
be in the case $d=n$ such that the dimension \is{of} domain and range
agree. But for the time being, we keep the possibility of having
different dimensions. In any case, we can use the standard basis vectors
$\ee_1,\ldots,\ee_n\in\R^n$ such that the above general theory reduces
to the following definition. 

\begin{definition}
  Given a probability space $(\Omega,\cA,\P)$\is{, a finite} number of
  vector-valued random variables 
$\gg_i = (g^{(i)}_1,\ldots, g^{(i)}_n)^\transpose:\Omega\to\R^n$, $1\le i\le N$,
  have a {\em joint Gaussian 
    distribution}  if the random variable $\GG:=(\gg_1,\ldots,
  \gg_N):\Omega\to \R^{n\times N}$ is uniquely determined by its matrix-valued mean
\[ 
 m^{(i)}_k = \E(g^{(i)}_k), \qquad 1\le i\le N, 1\le k\le n,
\]
and its tensor-valued covariance kernel
\[
\Sigma^{(ij)}_{k\ell} = 
\E[(g^{(i)}_k-m^{(i)}_k)(g^{(j)}_{\ell}-m^{(j)}_{\ell})], \qquad
  1\le i,j\le N,\quad 1\le k,\ell\le n.
\]
For $\cD\subseteq\R^d$, we will call a vector-valued random process
$\gg:\Omega\times\cD\to\R^n$ 
a {\em Gaussian process} if for every $N\in\N$ and every
$X=\{\xx_1,\ldots,\xx_N\}\subseteq \cD$, the vector-valued variables
$\gg_j=\gg(\cdot,\xx_j)$, $1\le j\le N$ have a joint Gaussian
distribution.  This is written as
$
\gg \sim {\rm GP}(\mm, \KK),
$
where the vector-valued function $\mm:\cD\to\R^n$ gives the mean function and
$K:\cD\times\cD\to \R^{n\times n}$ is the matrix-valued covariance
kernel, satisfying 
\[
\mm(\xx) = \E [\gg(\cdot,\xx)], \qquad
\KK(\xx,\xx') = \E[(\gg(\cdot,\xx)-\mm(\xx))(\gg(\cdot,\xx')-\mm(\xx'))^\transpose].
\]
Obviously, we have for fixed $i,j$ the relation $\KK(\xx_i, \xx_j) = 
[\Sigma^{(ij)}_{k\ell} : 1\le k,\ell\le n]$.
\end{definition}

We will also need a \tlg{vector-valued version} of the Karhunen-Lo\`eve
expansion. To this end, let $L_2(\cD)^n =
\{\ww = (w_1,\ldots, w_n) : w_k\in L_2(\cD)\}$ be the space of all
vector-valued square-integrable functions with inner product
$\langle \uu, \ww\rangle_{L_2(\cD)^n}:=\sum_{k=1}^n \langle
u_k,w_k\rangle_{L_2(\cD)}$. 

Noting that our covariance kernel $\KK:\cD\to\cD\to \R^{n\times n}$ of
the Gaussian process $\gg:\Omega\times\cD\to\R^n$ is 
not only continuous but also symmetric in the sense that $\KK(\xx,\xx')=
\KK(\xx',\xx)^\transpose$, we see that  the integral operator
$A:L_2(\cD)^n\to L_2(\cD)^n$, given by  
\begin{equation}\label{Aoperator-vector}
A\ww (\xx):=\int_\cD \KK(\xx,\yy) \ww(\yy) d\yy
\end{equation}
is a self-adjoint and compact operator. The operator is also positive,
as we have, assuming for simplicity a centered Gaussian process,
\begin{eqnarray*}
  \langle A\ww,\ww \rangle_{L_2(\cD)^n} & = & \sum_{k,\ell=1}^n \int_\cD
  \int_\cD \tlg{\KK}(\xx,\xx')_{k\ell} w_k(\xx) w_\ell(\xx') d\xx d\xx' \\
  & = & \sum_{k,\ell=1}^n \int_\cD
  \int_\cD \int_\Omega g_k(\omega,\xx) g_\ell(\omega,\xx') d\P(\omega)
  w_k(\xx) w_\ell(\xx') d\xx d\xx' \\
  & = & \int_\Omega \left|\int_\cD \sum_{k=1}^n g_k(\omega,\xx)
  w_k(\xx) d\xx\right|^2 d\P(\omega)\ge 0.
\end{eqnarray*}
Thus, the spectral theorem for compact, self-adjoint and positive
operators gives us again an orthonormal basis  of eigenfunctions of
$A$ with corresponding nonnegative eigenvalues. As before, we are only
interested in positive eigenvalues and \is{can} proceed as in the scalar-valued case to
derive Mercer's representation of the kernel and the Karhunen-Lo\`eve
expansion for vector-valued Gaussian processes.

\begin{proposition} \label{prop:KLV}
Let $\cD\subseteq\R^d$ be compact\is{, and let} $\gg:\Omega\times \cD\to\R^n$ be
an integrable, centered Gaussian process with continuous covariance
kernel $\KK:\cD\times\cD\to\R^{n\times n}$. Let $\{\pphi_k\}\subseteq
L_2(\cD)^d$ be an orthonormal subset of $L_2(\cD)^d$ consisting of
eigenfunctions of the integral operator associated to the kernel with
corresponding positive eigenvalues $\{\lambda_k\}$. Then, the kernel has the {\em
  Mercer representation}
\[
\KK(\xx,\xx') = \sum_{k=1}^\infty \lambda_k
\pphi_k(\xx)\pphi_k(\xx')^\transpose, \qquad \xx,\xx'\in \cD.
\]
The series converges uniformly and absolutely. Moreover, the Gaussian
\is{process} has the {\em Karhunen-Lo\`eve} representation
\begin{equation}\label{KLV-expansion}
\gg(\omega,\xx) = \sum_{k=1}^\infty \sqrt{\lambda_k} \xi_k(\omega)
\pphi_k(\xx), \qquad \xx\in\cD,\omega\in\Omega,
\end{equation}
where the convergence with respect to $\omega$ is in $L_2(\Omega,\P)$
and with respect to $\xx$ is uniform. The random variables
$\xi_k:\Omega\to\R$ are given \is{in terms of $g$} by
\[
\xi_k(\omega) = \frac{1}{\sqrt{\lambda_k}} \int_\cD \langle
\gg(\omega,\xx),\pphi_k(\xx)\rangle_2 d\xx.
  \]
\end{proposition}

With the notation above, we can now use vector-valued Gaussian
processes to form approximations to vector-valued functions
$\vv:\cD\subseteq\R^d\to\R^n$. As in the scalar-valued case we start
with a prior model $\vv_0:\Omega\times\cD\to\R^n$ with mean function
$\mm:\cD\to\R^n$, usually assumed to be zero, and a matrix-valued
covariance kernel $\KK:\cD\times\cD\to\R^{n\times n}$. Next, assuming
again that we have $N$ observations at design points
$X=\{\xx_1,\ldots,\xx_N\}\subseteq\cD$ in the form
$\yy_j=\vv(\xx_j)+\eepsilon_j$, $1\le j\le N$, where the $\eepsilon_j$
are either zero (exact data) or follow a normal distribution with mean
zero and variance $\sigma^2I_{n}$ with $I_{n}\in\R^{n\times n}$ being the identity
matrix and $\sigma^2>0$. In both cases, we can write the
conditional distribution \is{as a} posterior Gaussian process
$\vv_{N,\sigma}\sim {\rm GP}(\mm_{N,\sigma},\KK_{N,\sigma})$,
$\sigma\ge 0$. To describe the mean
$\mm_{N,\sigma}$ and the covariance $\KK_{N,\sigma}$ we will interpret $\mm(X):=
(\mm(\xx_1), \ldots, \mm(\xx_N)) \in \R^{n\times N} $as a vector in $\R^{nN}$
rather than a matrix in $\R^{n\times N}$. In the same fashion, we will
interpret the $N$ $(n\times n)$-matrices
$\KK(\xx,X) = (\KK(\xx,\xx_1), \ldots,
\KK(\xx, \xx_N))$ as a  matrix from $\R^{n\times nN}$ and the fourth
order tensor 
$\KK(X,X)=(\KK(\xx_i,\xx_j))_{1\le i,j\le N}$ consisting of $N^2$
  small $(n\times n)$ matrices as a matrix from $\R^{nN\times
    nN}$. Finally, we  set
  $\YY=(\yy_1^\transpose,\ldots,\yy_N^\transpose)^\transpose\in\R^{nN}$
  and use 
  $I\in\R^{nN\times nN}$ as the identity matrix.

  \begin{lemma}\label{lem:predictiveVectorCase}
   With the notation above, the predictive mean and covariance of the
   posterior Gaussian process  $\vv_{N,\sigma}$ conditioned by the
   data $\YY$ and $X$, has mean and covariance 
\begin{eqnarray*}
\mm_{N,\sigma}(\xx) &:=&\mm(\xx)  + \KK(\xx, X) (\KK(X,X)+\sigma^2
I)^{-1} \left(\YY- \mm(X)\right),\\ 
\KK_{N,\sigma}(\xx,\xx')&:=& \KK(\xx, \xx') - \KK(\xx,X)
(\KK(X,X)+\sigma^2I) ^{-1}\KK(\xx',X)^\transpose.
\end{eqnarray*}
In the case $\sigma=0$, we will again write $\mm_N$ and $\KK_N$, respectively.
  \end{lemma}
This is in accordance \is{with} scalar-valued Gaussian processes, and follows
immediately from the standard formulas of conditional Gaussian
processes.

\subsection{Divergence- and Rotation-free Gaussian Processes}

\is{Having introduced} the concept of vector-valued Gaussian processes, we
can now proceed to discuss divergence- and rotation-free Gaussian
processes. To this end, we recall that divergence and rotation of a
continuously differentiable  vector field $\vv:\cD\subseteq\R^d\to\R^d$ are defined by
\[
\div \vv = \nabla\cdot\vv= \sum_{k=1}^d \frac{\partial v_k}{\partial x_k}, \qquad
\]
for any $d\in\N$ and, for $d=2,3$, 
\[
\curl \vv = \nabla\times \vv = \begin{cases}
  \left(
  \frac{\partial v_3}{\partial x_2}-\frac{\partial v_2}{\partial
    x_3},
  \frac{\partial v_1}{\partial x_3}-\frac{\partial  v_3}{\partial
    x_1},
  \frac{\partial v_2}{\partial x_1}-\frac{\partial v_1}{\partial
    x_2}\right)^\transpose & \mbox{ for } d=3\\
  \frac{\partial v_2}{\partial x_1}-\frac{\partial v_1}{\partial
    x_2}\ & \mbox{ for } d=2.
\end{cases}
\]
\is{The} vector field is called divergence-free or rotation-free if
$\div \vv(\xx) = 0$ or $\curl \vv(\xx)=\00$, respectively. 

We can now
carry this concept over to vector-valued Gaussian processes.

\begin{definition}\label{def:divfree}
Let  $\gg:\Omega \times\cD\to \R^d$ be a Gaussian \is{process} defined on a compact
subset $\cD\subseteq\R^d$ which is continuously differentiable with
respect to the variable $\xx\in\R^d$.
\begin{itemize}
\item The process is called {\em divergence-free} if it satisfies
$ \div \gg(\cdot,\xx) = 0$ for all $\xx\in \cD$ and  $\P$-almost everywhere on $\Omega$.
\item The process is called  {\em rotation-free or curl-free} if 
$\curl \gg(\cdot,\xx) = 0$ is satisfied for all $\xx\in \cD$ and 0
$\P$-almost everywhere on $\Omega$.
\end{itemize}
In both cases, the differential operators act only with respect to the
spatial variable $\xx\in\cD$. 
\end{definition}

The property of being divergence-free or rotation-free carries over
from a Gaussian process to its mean and covariance kernel.

\begin{lemma} 
\label{lem:divfree}Let $\gg:\Omega\times \cD\to\R^d$ be a divergence-free Gaussian
process. Then the mean $\mm$ is also divergence-free, and the
matrix-valued covariance kernel $\KK:\cD\times\cD\to\R^{d\times d}$
has the following property. For a 
fixed $\xx'\in \cD$, each column of the map $\KK(\cdot,\xx')$ is in
$C^1(\cD)^d$ and divergence-free.

In the case of $d=2,3$ this statement holds \is{in a similar way} for a
rotation-free Gaussian process.
\end{lemma}

\begin{proof}
  We will prove the result only for a divergence-free Gaussian
  process.  We denote the components of $\gg$ by $g_k$, $1\le k\le
  d$. For the mean
\[
  \mm(\xx)= \E[\gg(\xx)] = \int_{\Omega}\gg(\omega,\xx) d\P(\omega),
\] 
we simply take the divergence under the integral, using the Leibniz theorem to justify
differentiating under the integral sign.  

For the covariance kernel, we first fix the column index $1\le \ell\le
d$. Then,  we have for the $k$-th component  of the $\ell$-th column 
 
\[
\tlg{\KK}(\xx, \xx')_{k \ell}
=  \int_\Omega
\left(g_{k}(\omega,\xx) - m_k(\xx)\right) \left(g_{\ell}(\omega,\xx') - m_\ell(\xx')\right)
d\P(\omega).   
\]

That the right-hand side \is{as a function of $\xx$ is}  in $C^1(\cD)$
for fixed $\xx'$ follows again from the Leibniz theorem. We can then
conclude that 
\begin{align*}
\div_{\xx} \KK(\xx,\xx')\ee_\ell 
&= \sum_{k=1}^d \frac{\partial}{\partial x_k} \tlg{\KK}(\xx, \xx')_{k\ell}\\
    & =  \sum_{k=1}^d \frac{\partial}{\partial x_k}\int_\Omega\!
\left(g_{k}(\omega,\xx) - m_k(\xx)\right) \left(g_{\ell}(\omega,\xx') - m_\ell(\xx')\right)
    d\P(\omega)\\
    & =  \int_\Omega \left[\div_{\xx} \left(\gg(\omega,\xx) -
      \mm(\xx)\right)\right] \left(g_{\ell}(\omega,\xx') -
    m_\ell(\xx')   \right)
    d\P(\omega)\\
& = 0,
\end{align*}
using the fact that both $\gg$ 
and $\mm$ are divergence-free.
\end{proof}

A Gaussian process is uniquely determined by its mean and covariance
kernel. This remains obviously true for vector-valued Gaussian
processes. As we have seen, a divergence-free Gaussian process has a
divergence-free mean and covariance kernel. Hence, it is natural to
ask the question whether it is enough for a Gaussian process to have a
divergence-free mean and covariance kernel to be divergence-free, as
well. The next theorem shows that this is indeed the case.

\begin{theorem}
  Let $\cD\subseteq\R^d$ be compact.
Let $\gg:\Omega\times\cD\to\R^d$ be a continuously differentiable vector-valued Gaussian
process whose mean is divergence-free and whose
covariance kernel has divergence-free columns.  Then
$  \gg$ is divergence-free.

An \is{analogous} statement holds \is{in the case $d=2,3$ for Gaussian processes with}
rotation-free mean and covariance.
\end{theorem}
\begin{proof}
  Again, we give the proof only for divergence-free case. Here, we
  will use  the vector-valued Karhunen-Lo\`eve expansion from
  Proposition \ref{prop:KLV} with $n=d$. We
  first have a look at the eigenfunctions $\pphi\in L_2(\cD)^d$ of the
  integral operator $A$. From (\ref{Aoperator-vector}) we have
  \[
  \lambda_k \pphi_k(\xx)= \int_{\cD} \KK(\xx,\yy) \pphi_k(\yy) d\yy.
  \]
\is{On taking the divergence of} each side and using the fact that the
divergence of each column of $\KK(\xx,\yy)$ with respect to $\xx$ is
zero \is{we see} that the eigenfunctions \is{$\phi_k$} are divergence-free.
Thus, the partial sums of the expansion (\ref{KLV-expansion}),
\[
 \sum_{k=1}^m \sqrt{\lambda_k} \xi_k(\omega)
\pphi_k(\xx),
\]
have divergence zero\is{, implying that the limit as $m \to \infty$ is also divergence-free}.
\end{proof}.

\subsection{Vector-Valued Reproducing Kernel Hilbert Spaces}

As in the scalar-valued case, there is a strong connection to
kernel-based approximation theory and reproducing kernel Hilbert
spaces. We start this discussion by introducing both concepts for \is{functions that are
Hilbert-space-valued}. After that, we look at the specific
cases of vector-valued functions and then at divergence- and rotation-free functions.

Let $W$ again be a Hilbert space and $\cL(W)$ the space of all bounded
operators from $W$ to $W$. As usual, for $A\in \cL(W)$ the {\em
  adjoint operator} $A^*\in \cL(W)$ is the unique operator satisfying
$\langle Av,w\rangle_W = \langle v, A^*w\rangle_W$ for all $v,w\in W$.

\begin{definition}
  \begin{enumerate}
  \item A continuous map $\KK:\cD\times\cD\to
    \cL(W)$ is called {\em positive definite} if it satisfies
    $\KK(\xx,\xx')^*=\KK(\xx',\xx)$ for all $\xx,\xx'\in\cD$ and 
    \[
    \sum_{i,j=1}^N \langle w_i, \KK(\xx_i,\xx_j)w_j \rangle_W >0.
    \]
    for all $N\in\N$\is{, all data sets $X=\{\xx_1,\ldots,\xx_N\}\subseteq\cD$ of
    pairwise distinct points,} and all $\ww\in W^N\setminus\{\00\}$.
    
  \item A Hilbert space $\cH:=\cH(\cD;W)=\{w:\cD\to W\}$ of $W$-valued
    functions is called a {\em reproducing kernel Hilbert space}
    (RKHS) if there is a function   $\KK:\cD\times \cD\to \cL(W)$  with
    \begin{enumerate}
    \item $\KK(\cdot,\xx) w \in \cH(\cD;W)$ for all $\xx\in \cD$ and
  all $w\in W$.
\item $\langle v(\xx),w\rangle_W = \langle v,
  \KK(\cdot,\xx)w\rangle_{\cH}$ for all $v\in\cH(\cD;W)$, all
  $\xx\in \cD$ and all $w\in W$.
\end{enumerate}
The function $\KK$ is called the {\em reproducing kernel} of $\cH(\cD;W)$.
  \end{enumerate}
\end{definition}

As in the case of Gaussian processes, this definition reduces to the
classical definition of positive definite kernels and reproducing
kernel Hilbert spaces for $W=\R$. While it is possible to draw
connections between the generalized concepts of Gaussian processes,
positive definite kernels and reproducing kernel Hilbert spaces, we
restrict ourselves here again to the case of vector-valued \is{processes}.

The restriction of the above definition to this case yields the
following concept of vector-valued RKHS and matrix-valued kernels.

\begin{definition}
  \begin{enumerate}
  \item A continuous, matrix-valued map $\KK:\cD\times\cD\to
    \R^{n\times n}$ is called {\em positive definite} if it satisfies
    $\KK(\xx,\xx')^\transpose=\KK(\xx',\xx)$ for all $\xx,\xx'\in\cD$ and 
    \begin{equation}\label{pdvectorkernel}
    \sum_{i,j=1}^N  \ww_i^\transpose \KK(\xx_i,\xx_j)\ww_j>0.
    \end{equation}
    for all $N\in\N$ and all data sets $X=\{\xx_1,\ldots,\xx_N\}\subseteq\cD$ of
    pairwise distinct points and all $\ww_1,\ldots,\ww_N\in\R^n$, not
    all being the zero vector.
    
  \item A Hilbert space $\cH:=\cH(\cD)=\{\ww:\cD\to\R^n\}$ of vector-valued
    functions is called a {\em reproducing kernel Hilbert space}
    (RKHS) if there is a function   $\KK:\cD\times \cD\to \R^{n\times n}$  with
    \begin{enumerate}
    \item $\KK(\cdot,\xx) \ww \in \cH(\cD)$ for all $\xx\in \cD$ and
  all $\ww\in \R^n$.
\item $\vv(\xx)^\transpose \ww = \langle \vv,
  \KK(\cdot,\xx)\ww\rangle_{\cH}$ for all $\vv\in\cH(\cD)$, all
  $\xx\in \cD$ and all $\ww\in \R^n$.
\end{enumerate}
The function $\KK$ is called the {\em reproducing kernel} of $\cH(\cD)$.
  \end{enumerate}
\end{definition}

As in the scalar-valued case, the reproducing kernel of a
vector-valued RKHS is unique and positive semi-definite, i.e. it
satisfies (\ref{pdvectorkernel}) with $\ge 0$ rather than $>0$. The
kernel is positive definite if point-evaluations are linearly
independent. Moreover, for every positive definite kernel
$\KK:\cD\times \cD\to\R^{n\times n}$ there is a unique vector-valued RKHS
$\cH=\cH(\cD)\subseteq C(\cD)^n$,  in which $\KK$ is the
reproducing kernel. It is, as in the scalar-valued case, constructed by
forming the closure of the linear space
\begin{equation}\label{FK}
F_{\KK}(\cD)  = \left\{\sum_{j=1}^M \KK(\cdot, \zz_j)\bbeta_j :
\zz_j\in \cD, \bbeta_j\in\R^n, M\in\N\right\}
\end{equation}
with respect to the norm which is induced by the inner product
\begin{equation}\label{FKinnerproduct}
\langle \KK(\cdot,\xx)\aalpha, \KK(\cdot,\xx')\bbeta\rangle_{\KK}:=
\aalpha^\transpose \KK(\xx,\xx')\bbeta.
\end{equation}
Details can, for example, be found in \cite{Fuselier-08-2}.

We will now describe the connection between vector-valued Gaussian
processes and vector-valued RKHS and matrix-valued positive definite
kernels. As in the scalar-valued case, we start with the no-noise situation,
meaning that our vector-valued observations at
$X=\{\xx_1,\ldots,\xx_N\}\subseteq\cD$ are given by  $\yy_j =
\vv(\xx_j)$, $1\le j\le N$ and thus we have $\sigma=0$ in Lemma
\ref{lem:predictiveVectorCase}.

\begin{theorem}\label{thm:matrixinterpolant}
  Let $\KK:\cD\times\cD\to\R^{n\times n} $ be a positive
  definite  kernel and  $X=\{\xx_1,\ldots,\xx_N\}\subseteq\cD$ consist of
  pairwise distinct points. Define the finite dimensional
  approximation space  $\tlg{V_X}=\{\sum_{j=1}^N \tlg{\KK}(\cdot,\xx_j)\bbeta_j :
  \bbeta_j\in\R^n, 1\le j\le N\}$ and  let $\cH=\cH(\cD)$ be the
  RKHS of the kernel $\KK$. 
  \begin{enumerate}
  \item For every set of training observations $\YY=\{\yy_1,\ldots,\yy_N\}\subseteq
    \R^n$ there is a unique function $\sb_{0,\YY}\in 
    \tlg{V_X}$ interpolating the data,
      i.e. satisfying $\sb_{0,\YY}(\xx_i) = \yy_i$, $1\le i\le N$. This defines
      an interpolation operator $\II_X:\R^{N\times n}\to \tlg{V_X}$, 
      $\II_X \vv=\sb_{0,\YY}$ and also an operator $\II_X:C(\cD)^n\to
      \tlg{V_X}$ by setting $\II_X\vv=\II_X(\vv|X)$. 
    \item The function $\sb_{0,\YY}\in V_X$ is the solution of
      \[
      \min\left\{\|\sb\|_\cH : \sb\in \cH \mbox{ with } \sb(\xx_i) = \yy_i,
      1\le i\le N\right\}.
      \]
    \item If $\mm_N$ is the predictive mean from Lemma
      \ref{lem:predictiveVectorCase} in the no-noise case
      $\yy_j=\vv(\xx_j)$, $1 \le j\le N$, i.e. with $\sigma=0$, then
      \begin{equation}\label{mnvector}
      \mm_N= \mm + \II_X \vv - \II_X \mm.
      \end{equation}
      \item If $\KK_N$ is the predictive covariance kernel from
        Lemma \ref{lem:predictiveVectorCase} in the no-noise case,
        i.e. with  $\sigma=0$, then
        \[ 
        \KK_N(\xx,\xx') = \KK(\xx,\xx') -\II_X \KK(\cdot,\xx') (\xx) =
        \KK(\xx,\xx') - \II_X \KK(\cdot,\xx) (\xx'),
        \]
        where the interpolation operator $\II_X$ is now applied to each
        column of the matrix-valued kernel. This again means that
        the predictive covariance kernel is the error between the
        given kernel and its interpolant with respect to one of its arguments.

        For the predictive variance
        $\ssigma_N^2(\xx):=\KK_N(\xx,\xx)\in\R^{n\times n}$
        this means
        \[
        \|\ssigma_N^2(\xx)\|_2= \lambda_{\max}(\ssigma_N^2(\xx))  =  \sup_{\|\vv\|_{\cH} =1} \|\vv(\xx) -\II_X\vv
        (\xx)\|_2^2,
        \]
        where $\lambda_{\max}(\ssigma_N^2(\xx))$ denotes the maximum
        eigenvalue of the positive and symmetric matrix
        $\KK_N(\xx,\xx)\in\R^{n\times n}$. 
  \end{enumerate}
\end{theorem}
\begin{proof}
  If we write the interpolant in the form
  \[
  \sb_{0,\YY} = \sum_{j=1}^N \KK(\cdot,\xx_j)\aalpha_j
  \]
  with $\aalpha_j\in\R^n$, then the interpolation conditions and the
  \is{notational convention immediately} before Lemma
  \ref{lem:predictiveVectorCase} show $\aalpha=\KK(X,X)^{-1}\YY$ and
  thus
  \[
  \sb_{0,\YY}(\xx) = \is{\KK(\xx,X)} \aalpha = \is{\KK(\xx,X)}
  \KK(X,X)^{-1}\YY,
  \]
proving the first point, from which also the third point follows
immediately. For the second point, we proceed similarly as 
in the scalar-valued case, see for example \cite{Wendland-05-1}. We
first note that if $\vv\in \cH$ is any function with $\vv(\xx_j)=\00$,
$1\le j\le N$, then
\[
\langle \vv, \sb_{0,\YY}\rangle_\cH = \sum_{j=1}^N\langle \vv,
\KK(\cdot,\xx_j)\aalpha_j\rangle_\cH
=\sum_{j=1}^N\vv(\xx_j)^\transpose\aalpha_j  = 0,
\]
using the reproducing kernel property. Thus, for any other interpolant
$\sb\in\cH$, $\sb(\xx_j)=\yy_j$, we can conclude
\[
\|\sb_{0,\YY}\|_\cH^2 = \langle \sb_{0,\YY}, \sb_{0,\YY}-\sb+\sb\rangle_{\cH}  = \langle
\sb_{0,\YY},\sb \rangle_\cH\le \|\sb_{0,\YY}\|_{\cH}\|\sb\|_{\cH}.
\]
Dividing by $\|\sb_{0,\YY}\|_{\cH}$ shows the second statement. For the
fourth statement, we first note that for $1\le j\le N$ and $1\le
\ell\le n$, there is a unique function $\uu_j^{(\ell)}\in \tlg{V_X}$
satisfying $\uu_j^{(\ell)}(\xx_i) = \delta_{ij}\ee_\ell$ for $1\le
i,j\le N$ and $1\le \ell\le n$, where $\ee_\ell$ is the $\ell$-th unit
vector in $\R^n$. Thus, \is{the matrix-valued function}
$\uu_j:=(\uu_j^{(1)},\ldots, \uu_j^{(n)}):\cD\to\R^{n\times n}$, \is{
satisfies} $\uu_j(\xx_i) = \delta_{ij} I_n$, where $I_n$ is the
identity matrix in $\R^{n\times n}$. The functions $\uu_j$ belong to
the finite dimensional space $\spn\{\KK(\cdot,\xx_j) : 1\le j\le
N\}$. Moreover, the functions $\uu_j^\transpose$ also satisfy
$\is{\uu_j(\xx_i)^\transpose} = \delta_{ij}I_n$ and belong to the space
$\spn\{\KK(\xx_j,\cdot) : 1\le j\le N\}$. Hence, they are, for a fixed
$\xx\in\cD$, the unique solution 
of the linear matrix system
\[
\sum_{j=1}^N \KK(\xx_i,\xx_j)\uu_j(\xx)^\transpose = \KK(\xx_i,\xx),
\qquad 1\le i\le N.
\]
With this, we can, again for a fixed $\xx\in\cD$ and $\aalpha\in\R^n$,
express the error as
\begin{eqnarray}
  \aalpha^\transpose [\vv(\xx)-\II_X\vv(\xx)] & = & \aalpha^\transpose
  \left[\vv(\xx) - \sum_{j=1}^N \uu_j(\xx) \vv(\xx_j)\right] \nonumber\\
  & = & \left\langle \vv, \left[\KK(\cdot,\xx)\aalpha - \sum_{j=1}^N
    \KK(\cdot,\xx_j)\uu_j(\xx)^\transpose
    \aalpha\right]\right\rangle_{\cH}\nonumber\\
  & \le & \|\vv\|_{\cH} P_{\KK,X}(\xx,\aalpha)\label{power1}
\end{eqnarray}
with
\begin{eqnarray*}
  [P_{\KK,X}(\xx,\aalpha)]^2  &=&   \left\|K(\cdot,\xx)\aalpha - \sum_{j=1}^N
    \KK(\cdot,\xx_j)\uu_j(\xx)^\transpose
    \aalpha\right\|_{\cH}^2\\
    & = & \aalpha^\transpose \KK(\xx,\xx)\aalpha - 2 \aalpha^\transpose\sum_{j=1}^N
    \KK(\xx,\xx_j)\uu_j(\xx)^\transpose \aalpha \\
&&\mbox{}    + \aalpha^\transpose \sum_{j,k=1}^N \uu_k(\xx)
    \KK(\xx_k,\xx_j)\uu_j(\xx)^\transpose \aalpha \\
    & = & \aalpha^\transpose \left[ \KK(\xx,\xx)-\sum_{j=1}^N
      \KK(\xx,\xx_j)\uu_j(\xx)^\transpose \right]\aalpha\\
    & = & \aalpha^\transpose \left[\KK(\xx,\xx)- \KK(\xx,X)
      \KK(X,X)^{-1}\KK(X,\xx)\right] \aalpha \\
    & = & \aalpha^T \KK_N(\xx,\xx)\aalpha.
\end{eqnarray*}

Obviously, for a fixed $\xx \in \cD$ we have equality in
(\ref{power1}) for
\[
\vv =
\KK(\cdot,\xx)\aalpha-\sum_{j=1}^N\KK(\cdot,\xx_j) \uu_j(\xx)^\transpose
\aalpha,
\] which implies
\[
\sup_{\substack{\vv\in\cH\\\vv\ne\00}}
\frac{|\aalpha^\transpose(\vv(\xx)-\II_X\vv(\xx))|^2}{\|\vv\|_{\cH}^2} =
P_{\KK,X}(\xx,\aalpha) ^2= \aalpha^\transpose \KK_N(\xx,\xx)\aalpha.
\]
This now leads to
\begin{eqnarray*}
  \lambda_{\max}(\KK_N(\xx,\xx)) & = &
\sup_{\substack{\aalpha\in\R^n\\\aalpha\ne \00}}
\frac{\aalpha^\transpose \KK_N(\xx,\xx)\aalpha}{\|\aalpha\|_2^2} \\
& = &
\sup_{\substack{\aalpha\in\R^n\\\aalpha\ne \00}}
\sup_{\substack{\vv\in\cH\\\vv\ne\00}}  \frac{|\langle 
  \vv(\xx)-\II_X\vv(\xx),\aalpha\rangle_2|^2}{\|\vv\|_{\cH}^2
  \|\aalpha\|_2^2}\\
& = &
\sup_{\substack{\vv\in\cH\\\vv\ne\00}}  
\sup_{\substack{\aalpha\in\R^n\\\aalpha\ne \00}} \frac{|\langle
  \vv(\xx)-\II_X\vv(\xx),\aalpha\rangle_2|^2}{\|\vv\|_{\cH}^2
  \|\aalpha\|_2^2}\\
& = &
  \sup_{\substack{\vv\in\cH\\\vv\ne\00}} \frac{ \|\vv(\xx)-\II_X\vv(\xx)\|_2^2}{\|\vv\|_{\cH}^2}
\end{eqnarray*}
\end{proof}

The connection between Gaussian processes and positive definite
kernels in the situation of noisy data is given by the following
vector-valued version of the representer theorem from machine
learning, Theorem \ref{thm:representer}.

\begin{theorem}\label{thm:matrixrepresenter}
  Let $\KK:\cD\times\cD\to\R^{n\times n} $ be a positive
  definite kernel and  $X=\{\xx_1,\ldots,\xx_N\}\subseteq\cD$ consist of
  pairwise distinct points. Define the finite dimensional
  approximation space  $\tlg{V_X}=\{\sum_{j=1}^N \KK(\cdot,\xx_j)\bbeta_j :
  \bbeta_j\in\R^n, 1\le j\le N\}$ and  let $\cH=\cH(\cD;\R^n)$ be the
  RKHS of the kernel $\KK$.
  \begin{enumerate}
    \item Given $\lambda>0$  and
      $Y=\{\yy_1,\ldots,\yy_N\}\subseteq\R^n$, there is a unique
      function $\sb_{\lambda,\YY}\in V_X$  solving
      \[
    \min\left\{\sum_{j=1}^N \|\sb(\xx_j)-\yy_j\|_2^2 + \lambda\|\sb\|_\cH^2 :
    \sb\in\cH\right\}.
    \]
  This defines an approximation operator $\QQ_{X,\lambda}:\R^{N\times
    n}\to \tlg{V_X}$ via
  $\QQ_{X,\lambda}\YY:=\sb_{\lambda,\YY}$ and also an operator
  $\QQ_{X,\lambda}:C(\cD)^n\to \tlg{V_X}$   
  by setting $\QQ_{X,\lambda}\vv := \QQ_{X,\lambda} (\vv|X)$.
  \item
  If $\mm_{N,\sigma}$ is the predictive mean from (\ref{eq:Gauss3}) with
  data $\YY$ and prior mean $\mm$ and covariance kernel $\KK$,
  then $\mm_{N,\sigma} = \mm-  \QQ_{X,\sigma^2}\mm +
  \QQ_{X,\sigma^2}(\YY)$. 
  \end{enumerate}
\end{theorem}
\begin{proof}
The proof of the first statement is with mild modification the same proof as the one for the
scalar-valued case. It is divided into two steps. First, we need to
show that the minimum is attained for an element $\sb\in \tlg{V_X}$. To see
this, we note that $\tlg{V_X}$ is a closed subspace of $\cH$ and hence we
can decompose $\cH$ into a direct sum of $\tlg{V_X}$ and its orthogonal
complement. Hence, any $\sb\in \cH$ can be written as
$\sb=\sb_1+\sb_2$ with $\sb_1\in \tlg{V_X}$ and $\sb_2\in \tlg{V_X}^\bot$. As we
have $\|\sb\|_{\cH}^2 = \|\sb_1\|_{\cH}^2+\|\sb_2\|_{\cH}^2$ and
$\aalpha^\transpose \sb_2(\xx_j) = \langle \sb_2,
\KK(\cdot,\xx_j)\aalpha\rangle_{\cH} = 0$ for all $\aalpha\in\R^n$ and
hence $\sb_2(\xx_j)=\00$ for $1\le j\le N$, we see that the minimum
indeed must be attained for an element from $\tlg{V_X}$. In the second step,
we write an arbitrary element of $\tlg{V_X}$ as
 $\sb=\sum_{j=1}^N \KK(\cdot,\xx_j)\bbeta_j$. Minimizing the above expression
then means determining  the
optimal $\bbeta_j$'s, which eventually gives the stated result.
\end{proof}

For us
the existence of  extension and restriction operators will be
important. If the kernel $\KK$ is actually defined on a larger set
$\widetilde{\cD}\supseteq \cD$, we can study the connection between
the reproducing kernel Hilbert space $\cH(\widetilde{\cD})$ and
$\cH(\cD)$ of $\KK$ and $\KK|\is{(\cD\times\cD)}$. The following result is a
generalization of the scalar-valued result, see \cite{Wendland-05-1}.

\begin{theorem}\label{thm:extension}
Let $\cD\subseteq\widetilde{\cD}\subseteq\R^d$ be given. Let
$\KK:\widetilde{\cD}\times \widetilde{\cD}\to\R^{n\times n}$ be positive
definite and let $\cH(\widetilde{\cD})$ and
$\cH(\cD)$ be the RKHS in which  $\KK$ and $\KK|\is{(\cD\times\cD)}$ are
reproducing kernels, respectively. 
\begin{itemize}
\item There is a linear isometric {\em extension operator} $E:\cH(\cD)\to
  \cH(\widetilde{\cD})$ with $\is{(E\vv)}|\cD=\vv$ for all $\vv\in \cH(\cD)$.
\item The linear {\em restriction operator}
  $R:\cH(\widetilde{\cD})\to \cH(\cD)$, defined by $R\vv=\vv|\cD$ for
  all $\vv\in \cH(\widetilde{\cD})$ satisfies $\|R\|=1$.
\end{itemize}
\end{theorem}

\begin{proof}
\is{Because} $\cD\subseteq\widetilde{\cD}$, we have a natural extension operator
$E: F_{\KK}(\cD)\to F_{\KK}(\widetilde{\cD})$ by simply evaluating
any $\vv\in F_{\KK}(\cD)$ also at points from $\widetilde{\cD}$. \is{(For the definition of $F_{\KK}$ see \eqref{FK}.)} This operator is obviously linear and an isometry,
\is{since} the norm of such a $\vv$ depends only on the centers and
coefficients of $\vv$. As $F_{\KK}(\cD)$ is dense in $\cH(\cD)$,
this extension operator extends to $\cH(\cD)$, keeping the isometry
property. 
  
We now turn to the restriction operator. Here, we need a vector-valued
version of the results of \is{Theorems 10.22, 10.26 and 
10.46} of \cite{Wendland-05-1}. We will adapt these results and
condense the arguments to the specific situation.
Hence, let $\vv\in \cH(\widetilde{\cD})$ be given. We first need to show
that $\vv|\cD\in\cH(\cD)$. To this end, we define a linear functional
$\psi_{\vv}:F_{\KK}(\cD)\to \R$ by setting $\psi_{\vv}(\KK(\cdot,\xx)\aalpha)
= \aalpha^\transpose \vv(\xx)$, $\xx\in\cD$, $\aalpha\in\R^n$, and
linearization. This is obviously linear. Moreover, for an arbitrary
$\ww=\sum_{j=1}^N \KK(\cdot,\xx_j)\aalpha_j\in F_{\KK}(\cD)$ we have
\begin{eqnarray*}
  \psi_{\vv}(\ww) & = & \sum_{j=1}^N \aalpha_j^\transpose \vv(\xx_j) =
  \sum_{j=1}^N \aalpha_j^\transpose \II_X\vv(\xx_j)
   =  \left\langle \II_X \vv, \sum_{j=1}^N
   \KK(\cdot,\xx_j)\aalpha_j\right\rangle_{\cH(\cD)}\\
   &=& \langle
  \II_X\vv,\ww\rangle_{\cH(\cD)}
   \le  \|\II_X\vv\|_{\cH(\cD)} \|\ww\|_{\cH(\cD)}
  \le \|\vv\|_{\cH(\widetilde{\cD})} \|\ww\|_{\cH(\cD)},
\end{eqnarray*}
where in the last step, we have used $\II_X\vv\in
F_{\KK}(\cD)\subseteq
F_{\KK}(\widetilde{\cD})\subseteq\cH(\widetilde{\cD})$ and the
norm-minimality of the interpolant. This shows boundedness of the
functional $\psi_{\vv}$ with $\|\psi_{\vv}\|\le
\|\vv\|_{\cH(\widetilde{\cD})}$. As $F_{\KK}(\cD)$ is dense in
$\cH(\cD)$, we can extend $\psi_{\vv}$ to a functional
$\psi_{\vv}:\cH(\cD)\to\R$ having the same norm. For simplicity, we will use the same
notation for the original functional and its extension. By the Riesz'
representation theorem, there is a $\ww_{\vv}\in \cH(\cD)$ with
$\|\ww_{\vv}\|_{\cH(\cD)}=\|\psi_{\vv}\|$ and $\psi_{\vv}(\ww) = \langle \ww,
\ww_{\vv}\rangle_{\cH(\cD)}$ for all $\ww\in \cH(\cD)$. It now
suffices to show that $\vv|\cD=\ww_{\vv}$. To see this, let
$\aalpha\in\R^n$ and $\xx\in\cD$ be fixed but arbitrary and conclude
\begin{eqnarray*}
  \aalpha^\transpose (\vv(\xx)-\ww_{\vv}(\xx)) & = &
  \aalpha^\transpose \vv(\xx) - \langle \ww_{\vv},
  K(\cdot,\xx)\aalpha)\rangle_{\cH(\cD)}\\
  &=& \aalpha^\transpose \vv(\xx) - \psi_{\vv}(K(\cdot,\xx)\aalpha) \\
  & = & \aalpha^\transpose \vv(\xx)-\aalpha^\transpose \vv(\xx) = 0.
\end{eqnarray*}
Now that we know that $R$ maps $\cH(\widetilde{\cD})$ indeed to
$\cH(\cD)$, we can look at the norm and see that for every $\vv\in
\cH(\widetilde{\cD})$ we have
\[
\|R\vv\|_{\cH(\cD)} = \|\ww_{\vv}\|_{\cH(\cD)} = \|\psi_{\vv}\| \le
\|\vv\|_{\cH(\widetilde{\cD})},
\]
showing $\|R\|\le 1$. However, any $\vv=\KK(\cdot,\xx)\aalpha$ with a
fixed $\xx\in \cD$ and $\aalpha\in\R^n$ belongs to $F_{\KK}(\cD)$
and
$F_{\KK}(\widetilde{\cD})\subseteq\cH(\widetilde{\cD})$. Obviously,
its restriction $R\vv$ is $\vv$ itself, showing $\|R\|=1$.
\end{proof}
  
\subsection{Sobolev Spaces}
In this paper, we are mainly interested in Sobolev spaces as
RKHS\tlg{s}. However, we will \is{use some leeway in} defining the inner product
and norm on these spaces, as \is{is customary} in the scalar-valued
case. The advantage of this procedure is that we can have different
reproducing kernels for the same algebraic space, using different
inner products. We start with recalling the definitions of
vector-valued Sobolev spaces and vector-valued mixed-regularity
Sobolev spaces. We do this not only for Hilbert spaces but \is{also for} general
$L_p$ spaces. To this end, we first recall the definition of the
$L_p$-norm for vector-valued functions. For $1\le p<\infty$ and
$\vv\in L_p(\cD)^n$, this is defined as
\[
\|\vv\|_{L_p(\cD)^n}:= \is{\left(\sum_{\ell=1}^n
\|v_{\ell}\|_{L_p(\cD)}^p\right)^{1/p} = \left(\sum_{\ell=1}^n \int_{\cD}
|v_{\tlg{\ell}}(\xx)|^p d\xx\right)^{1/p}},
\]
i.e. by taking the $\ell_p$ norm of the $L_p$ norms. For $p=\infty$
the definition of the norm \is{uses} the usual modifications.

For $k\in\N_0$, $1\le p<\infty$ the Sobolev space of order $k$ is
defined as expected:
 
\begin{eqnarray*}
  W_p^k(\cD)^n&:=&\{ \vv\in L_p(\cD)^n :
  \|\vv\|_{W_p^k(\cD)^n}<\infty\}, \\
  \|\vv\|_{W_p^k(\cD)^n} &:=& \is{\left(\sum_{\ell=0}^k |\vv|_{W_p^{\ell}(\cD)^n}^p\right)^{1/p}},\\
  \is{|\vv|_{W_p^{\ell}(\cD)^n}} &:=&   \is{\left(\sum_{|\aalpha|= \ell} \|D^\aalpha
  \vv\|_{L_p(\cD)^n}^2\right)^{1/2}}.
\end{eqnarray*}
As usual, the derivatives are meant in the weak sense. For $p=2$ we
will also use the notation $H^k(\cD)^n:=W_2^k(\cD)^n$. 
Moreover, for $1\le p<\infty$, using for example interpolation theory,
the definition can be extended to define also fractional order Sobolev
spaces $W_p^\tau(\Omega)^n$.

Before we specify the connection between Sobolev spaces and
reproducing kernel Hilbert spaces, we want to generalize a {\em
  sampling inequality} to the vector valued case. Such sampling
inequalities were introduced in the context of scattered data
interpolation and 
approximation in \cite{Narcowich-etal-05-1, Wendland-Rieger-05-1} and
then refined, for example, in
\cite{Arcangeli-etal-07-1,Arcangeli-etal-12-1,
  Narcowich-etal-06-1}. While an extension of 
these inequalities for functions vanishing on the data to
vector-valued functions is straight-forward and has been used, for
example, in \cite{Fuselier-08-1, Schraeder-Wendland-11-1,
  Wendland-09-1}, we require the following generalization for
functions not vanishing on the data. We formulate it only in the
required situation and not in the most general case. We will also need
a discrete $\ell_p$ norm, which we define for
$X=\{\xx_1,\ldots, \xx_N\}$ and $1\le p<\infty$ by
\[
\|\vv\|_{\ell_p(X)^n}:=\is{ \left(\sum_{\ell=1}^n
\|v_{\ell}\|_{\ell_p(X)}^p\right)^{1/p} = \left(\sum_{\ell=1}^n \sum_{j=1}^N
|v_{\ell}(\xx_j)|^p\right)^{1/p}}
\]
with the obvious modification for $p=\infty$. This, of course is a
generalization of the $\ell_p$ norm for vectors of vectors. For
$\YY=(\yy_1,\ldots,\yy_N)$ with $\yy_j\in\R^n$, we write
\[
\|\YY\|_{\ell_p}:= \left(\sum_{j=1}^{\is{N}}
\|\yy_j\|_{\ell_p}^p\right)^{1/p}.
\]

Recall that any bounded domain with a Lipschitz boundary also
satisfies an interior cone condition, see for example
\cite{Grisvard-85-1}. In the \is{sampling theorem that follows} we use the fill distance
from (\ref{filldistance}) and $(x)_+\is{:=}\max\{x,0\}$.

\begin{theorem}\label{thm:sampling}
Let $\cD\subseteq\R^d$ be a bounded domain with a Lipschitz
boundary. Let $p,q\in [1,\infty]$ \is{and $n \in \mathbb{N}$}. Let $\tau>d/2$ and
$\gamma:=\max\{2,p,q\}$. Then, there exist constants $h_0>0$
(depending on $\cD$ and $\tau$) and $C>0$ (depending on $\cD$,\tlg{n}, $\tau$ and $q$) such that for all $X=\{\xx_1,\ldots,\xx_N\}\subseteq \cD$
with $h=h_{X,\cD}\le h_0$ and all $\vv\in H^\tau(\cD)^n$,
\tlg{we have}
\[
|\vv|_{W_q^s(\cD)^n} \le C\left(h^{\tau-s-d(1/2-1/q)_+} |\vv|_{H^\tau(\cD)^n} +
h^{d/\gamma-s} \|\vv\|_{\ell_p(X)^n}\right)
\]
for all $0\le s\le \ell$, where $\ell=\ell_0:=\tau-d(1/2-1/q)_+$ in
the case $\tau\in\N$ and either $q>2$ and $\ell_0\in\N$ or
$q=2$. Otherwise $\ell=\lceil \ell_0\rceil-1$. Moreover, $s\in\N$ in the case
$q=\infty$.
\end{theorem}

\begin{proof}
In the case $n=1$, this is a special case of Theorems 3.1 and 3.2
from \cite{Arcangeli-etal-07-1}.  For $n\ge 2$ we set
 $e_1\tlg{:=}\tau-s-d(1/2-1/q)_+$ and  $e_2 \tlg{:=} d/\gamma-s$ to simplify the notation. Then, using the scalar-valued case, we have for $1\le
q<\infty$, 
\begin{eqnarray*}
  |\vv|_{W_q^s(\cD)^n} & = & \is{\left(\sum_{\ell=1}^n
  |v_{\ell}|_{W_q^s(\cD)}^q\right)^{1/q}}\\
  &\le & \is{\left(\sum_{\ell=1}^n \left[C\left(
    h^{e_1}|v_{\ell}|_{H^\tau(\cD)}+h^{e_2}\|v_{\ell}\|_{\ell_p(X)}\right)\right]^q\right)^{1/q}}\\
  &\le & \is{C\left[  h^{e_1} \left(\sum_{\ell=1}^n
  |v_{\ell}|_{H^\tau(\cD)}^q\right)^{1/q} + h^{e_2}\left(\sum_{\ell=1}^n \|v_{\ell}\|_{\ell_p(X)}^q\right)^{1/q}\right]}\\
  & \le& C\left[c_{2q} h^{e_1} |\vv|_{H^\tau(\cD)^n} + c_{pq} h^{e_2}
    \|\vv\|_{\ell_p(X)^n}\right],
\end{eqnarray*}
where we have used in the first step the scalar-valued bound, in the
second step the triangle inequality for $\ell_q$-norms and in the
final step the equivalence of $\ell_q$ and $\ell_p$ and $\ell_2$
norms. The equivalence constants are \is{given by} $c_{pq} = n^{1/q-1/p}$ if $p\ge q$ and $c_{pq}=1$ otherwise. The case $q=\infty$ is treated similarly.
\end{proof}

The Sobolev embedding theorem guarantees that for $\tau>d/2$, the
space $H^\tau(\cD)$  is a reproducing kernel Hilbert 
space. However, for general $\cD$ the reproducing kernel is usually
not known in explicit form. A typical remedy to this problem\is{, which we employ here,} is to use
a reproducing kernel of $H^\tau(\R^d)$, which can, for example, be
given as $K_\tau(\xx,\yy)=\Phi_\tau(\xx-\yy)$ with a continuous
and integrable function
$\Phi_\tau:\R^d\to\R$, which has a Fourier transform
\[
\widehat{\Phi}_\tau (\oomega) = (2\pi)^{-d/2}\int_{\R^d} \Phi_\tau(\oomega)
e^{-i\xx^\transpose \oomega} d\xx
\]
having the following decay property\is{: there exist} constants
$c_1,c_2>0$ such that
\begin{equation}\label{ft-decay}
c_1(1+\|\oomega\|_2^2)^{-\tau} \le \is{\widehat{\Phi}_\tau}(\oomega)\le
c_2(1+\|\oomega\|_2^2)^{-\tau},\qquad \oomega\in\R^d.
\end{equation}

To be more precise, we have the following well-known result see
\cite{Wendland-05-1}.

\begin{lemma}
Suppose $\Phi_\tau\in L_1(\R^d)\cap C(\R^d)$ possesses a Fourier
transform with decay (\ref{ft-decay}). Then, $K_\tau:\R^d\times
\R^d\to\R$, $K_\tau(\xx,\yy)=\Phi_\tau(\xx-\yy)$ is
the reproducing kernel of $H^\tau(\R^d)$ with respect to the inner
product
\[
\langle f,g\rangle_{K_\tau}:= (2\pi)^{-d/2} \int_{\R^d}
\frac{\widehat{f}(\oomega)
  \overline{\widehat{g}(\oomega)}}{\is{\widehat{\Phi}_{\tau}(\oomega)}} d\oomega.
\]
\end{lemma}

To use a kernel which is defined on all of $\R^d$ \is{and then} restricted to
$\cD\subseteq\R^d$ in the context of kernel-based learning can be done in two
equivalent ways. Both involve the existence of a bounded, linear
extension operator $E_\cD:H^\tau(\cD)\to H^\tau(\R^d)$, which
exists, for example, if $\cD$ is a bounded domain with a Lipschitz
boundary, see for example \cite{Brenner-Scott-08-1}. The first
approach is based on the fact that the interpolation/approximation
only uses data of the unknown function $f$ on data sites
$X\subseteq\cD$. Thus, the interpolation/approximation to $f$ is the
same as the interpolation/approximation to $E_\cD f$, which allows
us to work in the RKHS over all of $\R^d$. Here, we will use the
second approach, which uses the fact that there is also a canonical
extension operator for $E:\cH_K(\cD)\to \cH_K(\R^d)$ for any positive
definite kernel $K:\R^d\times \R^d\to\R$, which is the identity on
$F_K(\cD)=\spn\{K(\cdot,\xx) : \xx\in\cD\}$. This approach leads to
the following result; a proof can be found in \cite[Corollary
  10.48]{Wendland-05-1} for the integer case and follows by
interpolation theory for the fractional-order case.

\begin{lemma}\label{locns}
Assume $\Phi\in L_1(\R^d)\cap C(\R^d)$ has a Fourier
transform $\widehat{\Phi}$ satisfying (\ref{ft-decay}) with
$\tau>d/2$. Let $\cD\subseteq\R^d$ be a bounded domain with a
Lipschitz boundary. Let $K_\tau:\cD\times\cD\to\R$ be
defined by $K_\tau(\xx,\yy)=\Phi(\xx-\yy)$, $\xx,\yy\in\cD$. Then, there
exists an inner product $\langle
\cdot,\cdot\rangle_{K_\tau}:H^\tau(\cD)\times
H^\tau(\cD)\to\R$ on  $H^\tau(\cD)$ such that
$K_\tau$ is the reproducing
kernel of $H^\tau(\cD)$ with respect to this inner product. The
norm $\|\cdot\|_{K_\tau}$ induced by this inner product is
equivalent to the standard norm on $H^\tau(\cD)$, i.e. there
are constants $C_1(\tau),C_2(\tau)>0$ such that
\begin{equation}\label{Ctau}
C_1(\tau) \|u\|_{K_\tau} \le \|u\|_{H^\tau(\cD)} \le C_2(\tau)
\|u\|_{K_\tau}, \qquad u\in H^\tau(\cD).
\end{equation}
\end{lemma}
The inner product is first defined on $F_K(\cD):=\spn\{K(\cdot,\xx)
: \xx\in \cD\}$ by
\[
\langle K_\tau (\cdot,\xx), K_\tau
(\cdot,\yy)\rangle_{K_\tau}:=K_\tau(\xx,\yy)
\]
and linearisation and then extended to the closure of $F_K(\cD)$ with
respect to the norm induced by this inner product. For details, see \cite{Wendland-05-1}.
The above lemma can immediately be carried over to vector-valued
Sobolev spaces, using Theorem \ref{thm:extension}. However, we will do
tis for certain subspaces in the next section.

\subsection{Divergence- and Rotation-Free Kernels} 

In the context of kernel-based approximation theory, divergence- and
rotation-free kernels were  introduced in
\cite{Narcowich-Ward-94-2}. The construction of divergence-free or
rotation-free matrix-valued kernels is quite simple. Starting with a
positive definite  function $\Phi:\R^d\to\R$, which is at least twice
differentiable, we can define functions
$\PPhi_{\div}, \PPhi_{\curl}:\R^d\to\R^{d\times d}$ by 
\begin{eqnarray}
 \PPhi_{\div} &:=& (-\Delta I + \nabla\nabla^\transpose)\Phi,\label{phidiv}\\
 \PPhi_{\curl}&:=&-\nabla\nabla^\transpose \Phi.\label{phicurl} 
\end{eqnarray}
It is easy to see that for three times differentiable $\Phi$ the associated kernels
\is{
\begin{equation}\label{eq:KandPhi}
\KK_{\div}(\xx,\yy)=\PPhi_{\div}(\xx-\yy) \,\mbox{  and }\;\KK_{\curl}(\xx,\yy)=
\PPhi_{\curl}(\xx-\yy)
\end{equation}
}
have for a fixed \is{argument} $\yy$ divergence-free and
rotation-free columns respectively\is{, with the latter of course only for 
$d=2,3$}. 

These kernels are the reproducing kernels of function spaces consisting of
divergence- and rotation-free functions. To be more precise, the
following result comes from \cite{Fuselier-08-1, Fuselier-08-2,
  Schraeder-Wendland-11-1, Wendland-09-1}. To formulate it, we first
introduce the spaces of divergence- and rotation-free Sobolev
functions. For $\tau>1$ and $\cD\subseteq\R^d$, where $d=2,3$ in
the rotation-free case, let
\begin{eqnarray*}
  H^\tau(\cD;\div)&:=& \{\vv\in H^\tau(\cD)^d : \div \vv = 0 \}, \\
  H^\tau(\cD;\curl)&:=& \{\vv\in H^\tau(\cD)^d : \curl \vv = \00
  \}.
\end{eqnarray*}
Of course, as long as we do not have $\tau>d/2+1$, the derivatives are
meant to be in the weak sense and hence zero divergence and rotation
are meant only up to a set of measure zero. We will also need the
modified Sobolev spaces 
\begin{eqnarray*}
  \widetilde{H}^\tau(\R^d;\div) &=& \{\vv\in H^\tau(\R^d;\div) :
  \|\vv\|_{\widetilde{H}^\tau(\R^d)^d} < \infty\},\\
  \widetilde{H}^\tau(\R^d;\curl) &=& \{\vv\in H^\tau(\R^d;\curl) :
  \|\vv\|_{\widetilde{H}^\tau(\R^d)^d} < \infty\},
\end{eqnarray*}
which are both equipped with the inner product
 \begin{equation}\label{innerproduct1}
    \langle \vv,\ww\rangle_{\widetilde{H}^\tau(\R^d)^d}:=
    (2\pi)^{-d/2}\int_{\R^d} \frac{\overline{\widehat{\ww}(\oomega)}^\transpose
      {\widehat{\vv}(\oomega)}}{\|\oomega\|_2^2}
    (1+\|\oomega\|_2^2)^{\tau+1} 
    d\oomega.
 \end{equation}
Here, the Fourier transform of a vector-valued function is defined
component-wise. Note that for any $\vv\in
\widetilde{H}^\tau(\R^d;\div)$ or
$\vv\in\widetilde{H}^\tau(\R^d;\curl)$ we have the obvious but also
important bound $\|\vv\|_{H^\tau(\R^d)^d}
\le \|\vv\|_{\widetilde{H}^\tau(\R^d)^d}$.     

\tlg{
In this and the next section, we will use \eqref{ft-decay} with
$\tau$ replaced by $\tau+1$,
\begin{equation}\label{ft-decay-div0}
c_3(1+\|\oomega\|_2^2)^{-\tau-1} \le \widehat{\Phi}_{\tau+1}(\oomega)\le
c_4(1+\|\oomega\|_2^2)^{-\tau-1},\qquad \oomega\in\R^d.
\end{equation}
}
The connection between these modified Sobolev spaces and the above-defined kernels is as follows.
    
\begin{proposition}\label{prop:native1}
  Let $\tau>d/2$. Let $\Phi=\Phi_{\tau+1}$ define a reproducing kernel $K:\R^d\times\R^d\to\R$ of $H^{\tau+1}(\R^d)$ in the sense that \is{it}  has a Fourier transform $\widehat{\Phi}$  \tlg{satisfying (\ref{ft-decay-div0})}. 
  Assume that $\Phi\in W_1^2(\R^d)$. Then, 
  $H^{\tau+1}(\R^d)\subseteq C^1(\R^d)\cap W_\infty^1(\R^d)$ and $\Phi\in C^2(\R^d)$. \is{Moreover, the kernels $\KK_{\div}$ and $\KK_{\curl}$ defined from  $\PPhi_{\div}$ and $\PPhi_{\curl}$ as in  (\ref{eq:KandPhi})} are the
 reproducing kernels of $\widetilde{H}^\tau(\R^d;\div)$ and
 $\widetilde{H}^\tau(\R^d;\curl)$,
respectively, where in both cases the spaces are equipped with the inner product
    \[
    \langle \vv,\ww\rangle_{\widetilde{H}}=
    (2\pi)^{-d/2}\int_{\R^d} \frac{\overline{\widehat{\ww}(\oomega)}^\transpose
      {\widehat{\vv}(\oomega)}}{\widehat{\Phi}(\oomega)\|\oomega\|_2^2}
    d\oomega.
    \]
\end{proposition}

Next, we need extension operators to define kernel-induced inner
products on $H^\tau(\cD;\div)$ and $H^\tau(\cD,\curl)$. The following
result is from \cite{Wendland-09-1, Fuselier-08-1}.

\begin{proposition}\label{prop:extension}
 Let $\tau\ge 0$ and let $\cD\subseteq \R^d$ be a
 simply-connected, bounded domain with $d=2,3$.
 \begin{enumerate}
   \item If $\cD$ has   a $C^{\lceil\tau\rceil,1}$ 
  boundary then  there is a bounded,  linear operator
  $E_{\div}:H^\tau(\cD;\div)\to \widetilde{H}^\tau(\R^d;\div)$ with
  $(E_{\div}\vv)|\cD = \vv$ for all $\vv\in H^\tau(\cD;\div)$.
  \item If $\cD$ has a Lipschitz boundary then there exists a
    bounded and  linear operator $E_{\curl}:H^\tau(\cD;\curl)\to
    \widetilde{H}^\tau(\R^d;\curl)$ with $(E_{\curl}\vv)|\cD = \vv$
    for all $\vv\in H^\tau(\cD;\curl)$.
 \end{enumerate}
\end{proposition}

Note that we also need $d=2,3$ in the divergence-free case for the
existence of the extension operator. For that reason, we will, from
now on restrict ourselves to this situation.

With this, we can now prove a result \is{corresponding to} Lemma
\ref{locns} for our divergence-free and rotation-free kernels.

\begin{theorem} \label{thm:nativespace}
Let $d=2,3$. Assume $\Phi=\Phi_{\tau+1}\in W_1^2(\R^d)\cap C^2(\R^d)$ has a Fourier transform
  \tlg{satisfying (\ref{ft-decay-div0})} 
and $\tau>d/2$. Let  $\KK_{\div}$ be defined by  $\PPhi_{\div}$ from (\ref{phidiv}) and let
  $\KK_{\curl}$ be defined by $\PPhi_{\curl}$ from (\ref{phicurl}).
  Let $\cD\subseteq\R^d$ be a bounded, simply-connected domain.
  \begin{enumerate}
  \item If the domain $\cD$ has a $C^{\lceil\tau\rceil,1}$ boundary then there exists an
    inner product $\langle \cdot,\cdot\rangle_{\KK_{\div}}:
    H^\tau(\cD;\div)\times H^\tau(\cD;\div)\to \R$ such that
    $\KK_{\div}$ is the reproducing kernel with respect to this inner
    product. The norm defined by this inner product is equivalent to
    the standard norm on $H^\tau(\cD;\div)$.
  \item If the domain $\cD$ has a Lipschitz boundary then there is an
    inner product
    \(
    \langle \cdot,\cdot\rangle_{\KK_{\curl}}:
    H^\tau(\cD;\curl)\times H^\tau(\cD;\curl) \to\R
    \)
    such that
    $\KK_{\curl}$ is the reproducing kernel with respect to this inner
    product. The norm defined by this inner product is equivalent to
    the standard norm on $H^\tau(\cD;\curl)$.
  \end{enumerate}
\end{theorem}

\begin{proof}
As the proof of both cases is similar, we will only prove the first
one and write $\KK=\KK_{\div}$ to simplify the notation.  The proof
itself is similar to the proof of the scalar-valued case. We first
recall, that $\KK$ can be seen as the reproducing kernel of a unique
Hilbert space $\cH(\cD)$ but also as the reproducing kernel of a
Hilbert space $\cH(\R^d)$, both of which are constructed by completing the
pre-Hilbert spaces $F_{\KK}(\cD)$ and $F_{\KK}(\R^d)$, respectively. By Theorem
\ref{thm:extension}, there are natural extension and restriction
operators $E:\cH(\cD)\to \cH(\R^d)$ and $R:\cH(\R^d)\to\cH(\cD)$, both
having norm one. 

Next, on $\R^d$, we know by Proposition \ref{prop:native1} that the
spaces $\cH(\R^d)$ and $\widetilde{H}^\tau(\R^d;\div)$ are algebraically
the same with equivalent norms. More precisely, the identities
$\iota:{\cH(\R^d)\to \widetilde{H}^\tau(\R^d;\div)}$ and
$\widetilde{\iota}:{\widetilde{H}^\tau(\R^d;\div)\to \cH(\R^d)}$ satisfy
$\|\iota\|\le 1/C_1(\tau+1)^{1/2}$ and $\|\widetilde{\iota}\|\le
C_2(\tau+1)^{1/2}$, where the constants come from  (\ref{Ctau}).

Finally, we have the extension operator $E_{\div}:H^\tau(\cD;\div)\to
\widetilde{H}^\tau(\R^d;\div)$ from Proposition \ref{prop:extension}
and note that the restriction to $\cD$ of a function $\vv\in
\widetilde{H}^\tau(\R^d;\div)$  gives a function $\vv|\cD\in
H^\tau(\cD;\div)$ and thus defines a restriction operator
$R_{\div}:\widetilde{H}^\tau(\R^d;\div)\to H(\cD;\div)$ with
$\|\tlg{R_{\div}}\|=1$. 

Having all these operators at hand, we are now able to show that the spaces
$\cH(\cD)$ and $H^\tau(\cD;\div)$ are algebraically the same
spaces with equivalent norms. On the one hand, for $\vv\in
\cH(\cD)$ we have
\[
\vv = E\vv|\cD = \iota\circ E\vv|\cD = R_{\div}\circ\iota\circ E\vv
\in H^\tau(\cD;\div).
\]
On the other hand, for $\vv\in H^\tau(\cD;\div)$ we have
\[
\vv = E_{\div}\vv|\cD = \widetilde{\iota}\circ E_{\div}\vv |\cD =
R\circ\widetilde{\iota}\circ E_{\div} \vv \in \cH(\cD), \qquad
\]
showing that both spaces algebraically coincide. Moreover, these
relations also show
\begin{eqnarray*}
\|\vv\|_{H^\tau(\cD;\div)} &\le&
\frac{1}{C_1(\tau+1)^{1/2}} \|\vv\|_{\cH(\cD)},\\
\|\vv\|_{\cH(\cD)}&\le& C_2(\tau+1)^{1/2} C_{E}\|\vv\|_{H^\tau(\cD;\div)},
\end{eqnarray*}
where $C_{E}$ is the norm of $E_{\div}$. This is the stated norm
equivalence. 
\end{proof}

\section{Error Analysis}
\is{In this section we discuss the quality of the predictive mean as an estimate of the data-generating field $\vv: \cD \to \mathbb{R}^d$, under various assumptions on the noise.  Throughout we will make the following general assumptions on the geometry and the field $\vv$.

\begin{assumption}\label{ass1}
For $d = 2$ or $3$, $\cD \subseteq \mathbb{R}^d$ is a bounded domain with boundary satisfying 
\tlg{the conditions in Theorem \ref{thm:nativespace}, Part 1 or Part 2, as appropriate to the divergence-free or curl-free situations.}
\end{assumption}

\begin{assumption}\label{ass2}
 The field $\vv: \cD \to \mathbb{R}^d$ is a continuously differentiable function modelled by the Gaussian process \tlg{$\vv_0 \sim$}GP$(\mm, \KK)$, where $\KK$ is one of $\KK_{\div}$ or $
 \KK_{\curl}$ as defined in Theorem \ref{thm:nativespace}.
\end{assumption}

\begin{remark}
In practical situations it will often be desirable to modify the kernel  by the inclusion of ``\tlg{hyperparameters}'', for example replacing $\KK_{\rm div}(\xx, \xx')$ by 
$\alpha^2\KK_{\rm div}(\kappa\xx, \kappa\xx')$, where $\alpha$ and $\kappa$ are positive parameters, $\alpha^2$ controlling the variance and $\kappa$ the length scale.  For simplicity we ignore hyperparameters in the present discussion, noting that the predictive mean is independent of $\alpha^2$, while changing the inverse length scale parameter $\kappa$ would change only the constants in the error bounds in subsequent theorems.    
\end{remark}
}

\is{ \begin{assumption}\label{ass} The kernel $\KK$ 
\tlg{is the reproducing kernel of}
$H^\tau(\cD)^d$ and the target function $\vv$ \tlg{belongs} 
to $H^\beta(\cD)^d$\is, where $\beta, \tau > d/2$ in  order to ensure that functions in $H^\beta$ and $H^\tau $ are continuous. 
\end{assumption}
  
  Within this setting we  
   will consider both the {\em matching smoothness} case $\beta\ge \tau$
  and the {\em non-matching smoothness} case} $\tau>\beta$.

We will study the predictive mean computed under \is{both the assumption
  of no noise in the data, in which case we shoud look at $\mm_N$, and of noisy data, in which case we look  at $\mm_{N,\sigma}$. In the assumed no-noise case we
  will naturally derive error estimates when there is  no noise in the data, but we will
  also look at the misrepresented case in which we assume no noise but in fact there is noise in the data. In the case of assumed noisy 
  data \tlg{we} will
study first the misrepresented case that there is no noise, then the
  misrepresented case that there is noise but we do not know the kind
  of noise, and finally} the case that there is Gaussian noise.

For all of these cases, we will assume that the prior mean  $\mm$ is
divergence- or rotation-free, respectively, and  has the
same regularity as \is{  $\vv$. While the first assumption is
natural, a more natural smoothness assumption might} be that $\mm$ belongs to the RKHS of the
kernel used for the prior approximation, i.e. $\mm\in
H^\tau(\cD;\div)$, for example.  However, this choice \is{would} not
change the approximation results below significantly, as the
approximation order is always determined by $\min\{\beta,\tau\}$. \is{It
would only cause the statements of results to be lengthier}.

\subsection{Error Analysis for the Predictive Mean 
\tlg{using Interpolation}}

In this section, we discuss the error $\vv-\mm_N$, i.e. the error of
the predictive mean under the assumption that there is no noise in the
data, \tlg{even though there might be}\is{.  Thus} we will use the interpolation operator $\II_X$ to study the
error.

We start with the case that the observations $\yy_i=
\vv(\xx_i)$, $1\le i\le N$, indeed do not contain noise. We will state and
prove error estimates for \is{the} {\em matching case} that the
smoothness of $\vv$ is aligned with the smoothness of the kernel,
i.e. that $\vv\in H^\tau(\cD;\div)$ and $\tlg{\KK}_{\div}$ is the reproducing
kernel of the latter space in the sense of Theorem
\ref{thm:nativespace}. We will also have a result \is{in which} the
smoothness of the target function is less than required. In
approximation theory, this situation \is{has been}  named {\em escaping the
  native space}. As the only difference between the divergence-free
and the rotation-free case is the assumption on the smoothness of the
boundary of the domain, we will state both results together. 

We recall the fill distance from (\ref{filldistance})
and introduce the separation radius \is{and the mesh ratio}:
\begin{eqnarray*}
  h_{X,\cD} & = & \sup_{\xx\in\cD} \min_{\xx_j\in X} \|\xx-\xx_j\|_2,
  \\
  \is{q_{X,\cD}} & = & \frac{1}{2} \min_{\xx_j\ne \xx_k} \|\xx_j-\xx_k\|_2,\\ 
  \is{\rho_{X,\cD}} & = & \is{h_{X,\cD}/q_{X,\cD}}.
\end{eqnarray*}

To simplify the formulation of the approximation results, we will
employ the  notation
\begin{equation}\label{errorterm}
  \is{\mathcal{N}}_\is{\gamma}(\vv,\mm):=  \|\vv\|_{H^{\is{\gamma}}(\cD)^d} +
  \|\mm\|_{H^\is{\gamma}(\cD)^d}.
  \end{equation}
\begin{theorem}[Interpolation without
    Noise]\label{thm:interpolationnonoise} Let \is{Assumptions
    \ref{ass1},\ref{ass2} and \ref{ass} hold. Assume   \hw{that $\vv, \mm$
      belong to $H^\beta(\cD;\div)$} \tlg{or $H^\beta(\cD;\curl)$}.}
 Let $X=\{\xx_1,\ldots,\xx_N\}\subseteq\cD$ and
  let $\yy_j=\vv(\xx_j)$, $1\le j\le N$. Finally, let $\mm_N$ be the
  predictive, vector-valued mean based on the information
  $(\xx_j,\yy_j)$, $1\le j\le N$, \is{where the prior distribution uses} the
  kernel \is{$\KK = \KK_{\div}$ or $\KK_{\curl}$ and the} mean $\mm$. Then, there are constants
  $h_0,C>0$ such that for all  $X$ with $h_{X,\Omega}\le h_0$ and all $1\le q\le\infty$ 
    and all $0\le s\le \ell$ with $\ell$ given in Theorem
    \ref{thm:sampling} the following bounds hold.
  \begin{enumerate}
     \item
    If $\beta\ge \tau$ then 
    \[
    \|\vv-\mm_N\|_{W_q^s(\cD)^d} \le C h_{X,\cD}^{\tau-s -
      d(1/2-1/q)_+}\is{\mathcal{N}}_\tau(\vv,\mm).
    \]
      \item
    If $\beta<\tau$ then
    \[
    \|\vv-\mm_N\|_{W_q^s(\cD)^d} \le C h_{X,\cD}^{\beta-s -
      d(1/2-1/q)_+} \rho_{X,\cD}^{\tau-\beta} \is{\mathcal{N}}_\beta(\vv,\mm).
    \] 
  \end{enumerate}

\end{theorem}

\begin{proof}
  \is{For $\beta \ge \tau$} the statement is essentially \is{proved} as Theorem 5 in
  \cite{Fuselier-08-1}. However, we can shorten the proof by employing
  Theorem \ref{thm:nativespace}. We first note that with the interpolation
  operator $\II_X:C(\cD)^d\to \cH$ with $\cH=H^\tau(\cD;\div)$ or
  $\cH=H^\tau(\cD;\curl)$, respectively, from \is{(\ref{mnvector}) we have} 
  \begin{equation}\label{decomperror}
  \|\vv-\mm_N\|_{W_q^s(\cD)^d} \le \|\vv-\II_X\vv \|_{W_q^s(\cD)^d} +
  \|\mm-\II_X\mm\|_{W_q^s(\cD)^d}.
  \end{equation}
  As we assume that both $\vv$ and $\mm$ have Sobolev regularity
  $\beta\ge \tau>d/2$, both terms above are bounded in the same way
  and it suffices to concentrate on one of them, say the first one. As
  we have a sufficiently smooth boundary and as $\vv-\II_X\vv$
  vanishes on $X$, we can use \is{the sampling theorem, Theorem  \ref{thm:sampling},} to derive
  \[
  \|\vv- \II_X\vv\|_{W_q^s(\cD)^d} \le C h^{\tau-s-d(1/2-1/q)_+}
  \|\vv-\II_X\vv\|_{H^\tau(\cD)^d}.
  \]
  \is{By Theorem} \ref{thm:nativespace} we know  there is an
  equivalent norm on $\cH$ such that $\KK$ is the
  reproducing kernel of $\cH$ with respect to the
  corresponding inner product. Hence, by \is{the second part of Theorem
  \ref{thm:matrixinterpolant}, we have $\|\II_X v\|_{H^\tau(\cD)^d} \le C\|v\|_{H^\tau(\cD)^d}$ and hence $\|v - \II_X v\|_{H^\tau(\cD)^d} \le C \|\vv\|_{H^\tau(\cD)^d}$.} 

  For \is{$\beta < \tau$},  we still have (\ref{decomperror}) and the
  fact that $\vv$ and $\mm$ belong to the same Sobolev space. Hence,
  it again is enough to study one of the errors. Theorem
  \ref{thm:sampling} yields this time
   \[
  \|\vv- \II_X\vv\|_{W_q^s(\cD)^d} \le C h^{\beta-s-d(1/2-1/q)_+}
  \|\vv-\II_X\vv\|_{H^\beta(\cD)^d}.
  \]
  \is{We can} now proceed as in the proof of Theorem 6 in
  \cite{Fuselier-08-1} to use a specific auxiliary band-limited
  \is{approximation to $\vv$ together with the sampling theorem and a Bernstein inequality} to finally derive
  \begin{equation}\label{stability}
      \|\vv- \II_X\vv\|_{H^\beta(\cD)^d} \le C {\rho_{X,\cD}}^{\is{\tau-\beta}}
  \|\vv\|_{H^\beta(\cD)^d}.
  \end{equation}
\end{proof}
\begin{remark}
  Note that in the matching case, we do not need {\em quasi-uniform}
  point sets, i.e. point sets with $\rho_{X,\cD}\le \rho_0$,  as the
  error bound is solely expressed in terms of the 
  fill distance $h_{X,\cD}$. In the non-matching case, however,
  without quasi-uniformity, the ratio $\rho_{X,\cD}$ may diverge,
  leading to a \is{deterioration} of the approximation order.
  \end{remark}

Next, we want to discuss the interpolatory case, but now \is{with unexpected} noise in the data. Hence, instead of computing $\mm_N$ using
$\vv(\xx_j)$ we compute $\mm_N$ using the data
$\vv(\xx_j)+\eepsilon_j$. This leads to the following modified result.

\begin{theorem}[Interpolation with Noise]
  Under the \is{assumptions} of Theorem \ref{thm:interpolationnonoise}, let
  $\yy_j=\vv(\xx_j)+\eepsilon_j$, $1\le j\le N$\is{, and let $\mm_N$ be computed from }
  \tlg{Lemma~\ref{lem:predictiveVectorCase} with $\sigma=0$.}
  Then there are constants
  $h_0,C>0$ such that for all 
  $X$ with $h_{X,\Omega}\le h_0$ and all $1\le q\le\infty$ 
    and all $0\le s\le \ell$ with $\ell$ given in Theorem
    \ref{thm:sampling} the following bounds hold.
  
  \begin{enumerate}
     \item
    If $\beta\ge \tau$ then 
    \begin{eqnarray*}
 \lefteqn{\E\left[(\|\vv-\mm_N\|_{W_q^s(\cD)^d})\right] 
 \le C h_{X,\cD}^{\tau-s -
      d(1/2-1/q)_+}\is{\mathcal{N}} _\tau(\vv,\mm)}\\
&&\mbox{} 
+C\left(\tlg{h_{X,\cD}^{-s-d(1/2-1/q)_{+}+d/2}  
\rho^{\tau-d/2}_{X,\cD}} + h_{X,\cD}^{d/\gamma-s}\right)
  \is{\left(\E\left[\sum_{j=1}^N \|\eepsilon_j\|_2^2\right]\right)^{1/2}}.    
\end{eqnarray*}

\item
If $\beta<\tau$ then
\begin{eqnarray*}
&&    \E\left[(\|\vv-\mm_N\|_{W_q^s(\cD)^d})\right]
      \le C h_{X,\cD}^{\beta-s -
      d(1/2-1/q)_+} \rho_{X,\cD}^{\tau-\beta} \is{\mathcal{N}}_\beta(\vv,\mm)\\
&&\mbox{}  
+C\left(\tlg{h_{X,\cD}^{ -s-d(1/2-1/q)_{+}+d/2}  
\rho^{\tau-d/2}_{X,\cD}} + h_{X,\cD}^{d/\gamma-s}\right)
    \is{\left(\E\left[\sum_{j=1}^N
\|\eepsilon_j\|_2^2\right]\right)^{1/2}}.
\end{eqnarray*}
\is{where $\gamma=\max\{2,q\}$, and  } the expectation is taken with
respect to the distribution of the errors $\eepsilon_j$.
  \end{enumerate}
\end{theorem}

\begin{proof}
In this situation, we have $\mm_N=\mm-\II_X\mm + \II_X\vv + \II_X\cEE$ with
$\cEE=(\eepsilon_1,\ldots,\eepsilon_N)$. Thus, the bound
(\ref{decomperror}) on  $\|\vv-\mm_N\|_{W_q^s(\cD)^d}$ now becomes
\[
\|\vv-\mm_N\|_{W_q^s(\cD)^d} \le \|\vv-\II_X\vv \|_{W_q^s(\cD)^d} +
\|\mm-\II_X\mm\|_{W_q^s(\cD)^d} + \|\II_X\cEE\|_{W_q^s(\cD)^d}.
\]
As we have already bounded the first two terms, we only need to derive
a bound on the interpolant of the error. As the interpolant belongs to
$H^\tau(\cD)^d$, we can use the sampling inequality in the form
\begin{equation}\label{ineq_3rd_term}
\|\II_X\cEE\|_{\is{W_q^s(\cD)^d}} \le\! C\left(
h^{\tau-s-d(1/2-1/q)_+}\|\II_X\cEE\|_\is{{H^\tau(\cD)^d}} + \is{h_{X,\cD}}^{d/\gamma-s}
  \|\II_X\cEE\|_{\ell_2(X)^d}\right).
\end{equation}
On the one hand, we obviously have
$\|\II_X\cEE\|_{\ell_2(X)^d}^2=\sum_{j=1}^{\tlg{N}} \|\eepsilon_j\|_2^2$. On the
other hand, using again the equivalence of norms from Theorem
\ref{thm:nativespace}, we see \is{that}
\[
\|\II_X\cEE\|_{H^\tau(\cD)}^2 \le C \cEE^\transpose \KK(X,X)^{-1}\cEE \le C
q_X^{-2\tau+d}\|\cEE\|_2^2,
\]
where the last inequality follows from \cite[Theorems 6 and
  7]{Fuselier-08-2}, which state that there is a constant $C>0$ such
that the smallest eigenvalue of the interpolation matrix can be
bounded as
\[
\lambda_{\min}(\KK(X,X)) \ge C q_X^{-d-2}\inf_{\|\oomega\|_2\le 1/q_X}
\widehat{\Phi}(\oomega).
\]

\tlg{
Hence, inequality \eqref{ineq_3rd_term} can be bounded by
\begin{align*}
\|\II_X\cEE\|_{W_q^s(\cD)^d}&\le C
h^{\tau-s-d(1/2-1/q)_+}_{X,\cD} q_X^{-\tau+d/2}\|\cEE\|_2  + Ch_{X,\cD}^{d/\gamma-s}
  \left(\sum_{j=1}^N \|\eepsilon_j\|_2^2\right)^{1/2}\\
&\le C(h_{X,\cD}^{-s-d(1/2-1/q)_{+}+d/2}  
\rho^{\tau-d/2}_{X,\cD} + h_{X,\cD}^{d/\gamma-s})
  \left(\sum_{j=1}^N \|\eepsilon_j\|_2^2\right)^{1/2}.
\end{align*}
We obtain the result by combining the estimates from Theorem~\ref{thm:interpolationnonoise} and the above estimate.
}
\end{proof}

\subsection{Error Analysis for the Predictive Mean using Approximation}

In this section, we will assume that the measurements at
$X=\{\xx_1,\ldots,\xx_N\}\subseteq\cD$ of a vector-valued function
$\vv:\cD\to\R$ are contaminated, i.e. that we observe
$\yy_j=\vv(\xx_j)+\eepsilon_j$, $1\le j\le N$.
However, for the time being, we will not make any assumption on the
distribution of the noise, even allowing no noise at all. The
difference to the previous section is that we \is{now use not the interpolation operator but rather} the operator
$\is{\QQ_{X,\lambda}}:C(\cD)\to H^\tau(\Omega;\div)$ from Theorem
\ref{thm:matrixrepresenter} with a general smoothing parameter
$\lambda>0$ and the kernel $\KK_{\div}$. The
approximation operator indeed maps to $H^\tau(\Omega;\div)$ by Theorem
\ref{thm:nativespace}. Of course, the same holds true in the
rotation-free case.

To apply the sampling inequality from Theorem
\ref{thm:sampling}, we need the following bounds which are the
vector-valued version of similar bounds in the scalar-valued case.

\begin{lemma}\label{lem:standardbounds}
  Let the assumptions on the domain and kernel of Theorem
  \ref{thm:nativespace} hold. Let
  $\KK=\KK_{\div}$ or $\KK=\KK_{\curl}$, and $\cH=H^\tau(\cD,\div)$ or
  $\cH=H^\tau(\cD,\curl)$, respectively, equipped with the inner
  product in which $\KK$ is the reproducing kernel. Let
  $X=\{\xx_1,\ldots,\xx_N\}\subseteq{\cD}$ and
  $\YY=(\yy_1,\ldots,\yy_N)$ be given by 
  $\yy_j=\vv(\xx_j)+\eepsilon_j$, $1\le j\le N$, with a function
  $\vv\in\cH$ and including the possibility $\eepsilon_j=\00$. If
  $\QQ_{X,\lambda}$ is the approximation operator from Theorem
  \ref{thm:matrixrepresenter} and
  $\cEE=(\eepsilon_1,\ldots,\eepsilon_N)$ then
  \begin{eqnarray*}
    \|\YY-(\QQ_{X,\lambda}\YY)|X\|_{\ell_2} & \le & 
    \left(\|\cEE\|_{\ell_2}^2+\lambda \|\vv\|_{\cH}^2\right)^{1/2},\\
\|\QQ_{X,\lambda} \YY\|_{\cH} & \le &
\left(\lambda^{-1}\|\cEE\|_{\ell_2}^2 + \|\vv\|_{\cH}^2\right)^{1/2}.
\end{eqnarray*}
\end{lemma}
\begin{proof}
As $\QQ_{X,\lambda}\YY$ minimizes the penalized least-squares expression,  we  have
\begin{eqnarray*}
  \lefteqn{
  \is{
    \|\YY-(\QQ_{X,\lambda}\YY)|X\|_{\ell_2}^2 + \lambda\|\QQ_{X,\lambda}\YY\|_{\cH}^2}}
    \\
& \is{=} & \is{\sum_{j=1}^N   \|\yy_j-\QQ_{X,\lambda}\YY(\xx_j)\|^2_\tlg{2}}  +
\lambda\|\QQ_{X,\lambda}\YY\|_{\cH}^2 \\
& \le & \is{\sum_{j=1}^N   \|\yy_j-\vv(\xx_j)\|^2_\tlg{2}}  +
\lambda\|\vv\|_{\cH}^2\\
& \is{=} & \is{\sum_{j=1}^N \|\eepsilon_j\|^2_\tlg{2}} + \lambda \|\vv\|_{\cH}^2\\
& = & \|\cEE\|_{\ell_2}^2 + \lambda\|\vv\|_{\cH}^2,
\end{eqnarray*}
which proves both statements.
\end{proof}

We will also need the following auxiliary result for considering
approximation with mismatched smoothness. To state it, we need to
recall {\em band-limited} functions, which are functions having a
compactly supported Fourier transform. For a parameter $\delta>0$, we
will use the spaces
\begin{eqnarray*}
  \cB^\delta &:=& \{\vv\in L_2(\R^d)^d : \supp \tlg{\widehat{\vv}} \subseteq
  B_\delta[0]\},\\
    \widetilde{B}^\delta &:=& \left\{\vv\in\cB_\delta : \int_{\R^d}
  \frac{\|\widehat{\vv}(\oomega)\|_2^2}{\|\oomega\|_2^2} d\oomega
  <\infty\right\},
  \end{eqnarray*}
where $B_\delta[0]$ denotes the closed unit ball \is{in $\R^d$} with radius
$\delta>0$. Moreover, we need the sub-spaces 
$\widetilde{\cB}^\delta_{\div}$ and $\widetilde{\cB}^\delta_{\curl}$
of $\widetilde{\cB}^\delta$, containing the divergence-free and
rotation-free vector fields, respectively. The following lemma provides a
simple scaling argument as the first result and then 
summarizes Theorems 1 and 2 of \cite{Fuselier-08-1}, noting that $t=0$
 is also possible there, \is{as was} the case in \cite[Theorem
  3.4]{Narcowich-etal-06-1}. 

\begin{lemma}\label{lem:bandlimited}
\begin{enumerate}
\item If $\tau\ge \beta\ge 0$, $\delta\ge 1$ and $\vv\in
  \widetilde{\cB}^\delta_{\div}$ then
  \[
  \|\vv\|_{\widetilde{H}^\tau(\R^d)^d} \le
  2^{(\tau-\beta)/2}\delta^{\tau-\beta}\|\vv\|_{\widetilde{H}^\beta(\R^d)^d}.
  \]
\item Let $\tau\ge\beta>d/2$. Given $\vv\in \widetilde{H}^\tau(\R^d;\div)$
and a point set $X\subseteq\R^d$  with separation radius $q_X$, there
exists a function $\vv_\delta\in \widetilde{\cB}^\delta_{\div}$ with
$\vv_{\delta}|X=\vv|X$ and
\[
\|\vv-\vv_\delta\|_{\widetilde{H}^\beta(\R^d)^d} \le
5\kappa^{\beta-\tau} q_X^{\tau-\beta}\|\vv\|_{\widetilde{H}^\tau(\R^d)^d},
\]
with $\delta=\kappa/q_X$, where $\kappa\ge 1$ depends only on $\beta$
and $d$.
\item The same statements hold for rotation-free functions.
\end{enumerate}
\end{lemma}
  
We start our investigation on the approximation error by first falsely
assuming that there is an \is{error} though the values are actually
exact. This leads to the following result. As before, we assume that
the prior mean has the same regularity as the unknown function.

\begin{theorem}[Approximation without
    Noise]\label{thm:approximationnonoise}Let \is{Assumptions
    \ref{ass1}, \ref{ass2} and \ref{ass} hold. Assume   \hw{$\vv, \mm$
      belong to $H^\beta(\cD;\div)$} \tlg{or $H^\beta(\cD;\curl)$}.}
 Let $X=\{\xx_1,\ldots,\xx_N\}\subseteq\cD$ and
  let $\yy_j=\vv(\xx_j)$, $1\le j\le N$.Finally, \is{given $\lambda > 0$} let $\mm_{N,\lambda}$ be the
  predictive, vector-valued mean, based on the information
  $(\xx_j,\yy_j)$, $1\le j\le N$, and the prior distribution using the
  kernel \is{$\KK = \KK_{\div}$ or $\tlg{\KK_{\curl}}$} and the initial mean $\mm$. Then, there are constants
  $h_0,C>0$ such that for all  $X$ with $h_{X,\Omega}\le h_0$, all $1\le q\le\infty$ 
    and all $0\le s\le \ell$ with $\ell$ given in Theorem
    \ref{thm:sampling} the following bounds hold.
  \begin{enumerate}
     \item
    If $\beta\ge \tau$ then 
\[
    \|\vv-\mm_{N,\lambda}\|_{W_q^s(\cD)^d} \le C\left( h_{X,\cD}^{\tau-s -
      d(1/2-1/q)_+} + \tlg{\sqrt{\lambda}} h_{X,\cD}^{d/\gamma-s}\right)\is{\mathcal{N}}_\tau(\vv,\mm).
\]
      \item
    If $\beta<\tau$ then
    \begin{eqnarray*}
    \|\vv-\is{\mm_{N, \lambda}}\|_{W_q^s(\cD)^d} &\le& C h_{X,\cD}^{\beta-s -
      d(1/2-1/q)_+} \rho_{X,\cD}^{\tau-\beta} \is{\mathcal{N}}_\beta(\vv,\mm)\\
    &&\mbox{}
     + C \tlg{\sqrt{\lambda}} h_{X,\cD}^{\beta-\tau +
      d/\gamma-s}\rho_{X,\cD}^{\tau-\beta} \is{\mathcal{N}}_\beta(\vv,\mm),
    \end{eqnarray*}
\is{where} $\gamma = \max\{2,q\}$.
  \end{enumerate}

\end{theorem}
\begin{proof}
As we assume no noise, we have
$\mm_{N,\lambda}=\mm-\is{\QQ_{X,\lambda}}\mm + \is{\QQ_{X,\lambda}}\vv$,
which allows us to split the error again as
\[
\|\vv-\mm_{N,\lambda}\|_{W_q^s(\cD)^d} \le
\|\vv-\is{\QQ_{X,\lambda}}\vv\|_{W_q^s(\cD)^d} +
\|\mm-\is{\QQ_{X,\lambda}}\mm\|_{W_q^s(\cD)^d}.
\]
As we assume $\vv$ and $\mm$ to have the same regularity, we
concentrate again on $\vv$. Here, the sampling inequality from Theorem
\ref{thm:sampling} and Lemma \ref{lem:standardbounds} with $\cEE=\00$ immediately
yield in the case $\beta\ge \tau$, 
\begin{eqnarray*}
  \|\vv-\is{\QQ_{X,\lambda}}\vv\|_{W_q^s(\cD)^d} & \le &
  C h^{\tau-s-d(1/2-1/q)_+}\|\vv-\is{\QQ_{X,\lambda}}\vv
  \|_{H^\tau(\cD)^d}\\
  & & \mbox{} +
  C h^{d/\gamma-s}\|\vv-\is{\QQ_{X,\lambda}}\vv\|_{\ell_2(X)^d}\\
  & \le & C h^{\tau-s-d(1/2-1/q)_+}\|\vv \|_{H^\tau(\cD)^d}\\
   & & \mbox{} +
  C h^{d/\gamma-s}\tlg{\sqrt{\lambda}} \|\vv\|_{H^\tau(\cD)^d},
\end{eqnarray*}
where we have set \is{$h=h_{X,\Omega}$}.

For the second statement, we can again concentrate on $\vv$, as the
same arguments apply to $\mm$.  For the extended function
$E_{\div}\vv\in \widetilde{H}^\beta(\R^d;\div)$, we choose an interpolatory
band-limited function $\vv_\delta\in  
 \widetilde{B}^\delta_{\div}$ with $\delta=\kappa/q_X$, according to
 the second point of Lemma \ref{lem:bandlimited}. Then, we start with
\begin{equation}\label{split2}
\|\vv-\is{\QQ_{X,\lambda}}\vv\|_{W_q^s(\cD)^d}\le
\|\vv-\vv_\delta\|_{W_q^s(\cD)^d}+
\|\vv_\delta-\is{\QQ_{X,\lambda}}\vv\|_{W_q^s(\cD)^d}.
\end{equation}
To bound the first term,  we use the  sampling inequality from Theorem
\ref{thm:sampling}, noting $\vv|X=\vv_\delta|X$, yielding
\begin{eqnarray*}
 \|\vv- \vv_\delta\|_{W_q^s(\cD)^d} &\le& C h^{\beta-s-d(1/2-1/q)_+}
 \|\vv-\vv_\delta\|_{H^\beta(\cD)^d}\\
 &\le& C  h^{\beta-s-d(1/2-1/q)_+}\|\vv\|_{H^\beta(\cD)^d},
 \end{eqnarray*}
where  the last inequality follows via 
\[
\|\vv- \vv_\delta\|_{H^\beta(\cD)^d} \le \|E_{\div}\vv -
\vv_\delta\|_{\widetilde{H}^\beta(\R^d)^d} \le \is{C}\|E_{\div}\vv\|_{\widetilde{H}^\beta(\R^d)^d}
\le C \|\vv\|_{H^\beta(\cD)^d},
\]
where we have used
Lemma \ref{lem:bandlimited}
with $\beta=\tau$ and the boundedness of the extension operator.

For the second term in (\ref{split2}), we 
  first note that we have $\is{\QQ_{X,\lambda}}\vv
 = \is{\QQ_{X,\lambda}}\vv_\delta$, as the band-limited function has the same values on
 $X$ as the function and the approximation operator uses only these
 values. Thus, we can once again use the sampling inequality to derive
 \begin{eqnarray*}
 \|\vv_\delta-\is{\QQ_{X,\lambda}}\vv_\delta\|_{W_q^s(\cD)^d}&\le&
 h^{\tau-s-d(1/2-1/q)_+}
 \|\vv_\delta-\is{\QQ_{X,\lambda}}\vv_\delta\|_{H^\tau(\cD)^d}\\
 & & \mbox{} +
     C h^{d/\gamma-s}\|\vv_\delta-\is{\QQ_{X,\lambda}}\vv_\delta\|_{\ell_2(X)^d}.
 \end{eqnarray*}
To further bound the norms on the right-hand side, we use Lemma
\ref{lem:standardbounds} in combination with Lemma \ref{lem:bandlimited} and the
equivalence of norms from Theorem \ref{thm:nativespace}. This yields on the one hand
\begin{align*}
\|\vv_\delta-\is{\QQ_{X,\lambda}}\vv_\delta\|_{H^\tau(\cD)^d} &\le
\is{ C\|\vv_\delta\|_{H^\tau(\cD)^d}}
\le \is{C}\|\vv_{\delta}\|_{\widetilde{H}^\tau(\R^d)^d}
\le C
q_X^{\beta-\tau}\|\vv_\delta\|_{\widetilde{H}^\beta(\R^d)^d}\\
&\le C h^{\beta-\tau}
\rho^{\tau-\beta} \|\vv\|_{H^\beta(\cD)^d},
\end{align*}
where the last inequality follows from
\begin{eqnarray*}
\|\vv_\delta\|_{\widetilde{H}^\beta(\R^d)^d} &\le& \|\vv_\delta-
E_{\div}\vv\|_{\widetilde{H}^\beta(\R^d)^d} +
\|E_{\div}\vv\|_{\widetilde{H}^\beta(\R^d)^d}
\le C \|E_{\div}\vv\|_{\widetilde{H}^\beta(\R^d)^d}\\
&\le& C \|\vv\|_{H^\beta(\cD)^d}.
\end{eqnarray*}
On the other hand, \is{from Lemma \ref{lem:standardbounds} with $\cE = 0$} we have
\[
\|\vv_\delta - \is{\QQ_{X,\lambda}}\vv_\delta\|_{\ell_2(X)^d} \le
C \tlg{\sqrt{\lambda}} \|\vv_\delta\|_{H^\tau(\cD)} \le C \tlg{\sqrt{\lambda}}
h^{\beta-\tau}\rho^{\tau-\beta} \|\vv\|_{H^\beta(\cD)}.
\]

\end{proof}

\begin{remark} The above result shows that we can choose the {\em
    smoothing parameter} $\lambda$ as a function of $h_{X,\cD}$ to
  obtain the same convergence results as in the interpolatory case.
\end{remark}

Next, we discuss the approximation properties of the approximation
operator in the case of noisy data in the most general situation,
i.e. we do not make any assumptions on the noise.

\begin{theorem}[Approximation with
    Noise]\label{thm:approximationnoise} Let the assumptions of
  Theorem \ref{thm:approximationnonoise} hold, except that now the data are
  given by $\yy_j=\vv(\xx_j)+\eepsilon_j$, $1\le j\le N$. Again, \is{given $\lambda > 0$,} let $\mm_{N,\lambda}$ be the
  predictive, vector-valued mean, based on the information
  $(\xx_j,\yy_j)$, $1\le j\le N$, and the prior distribution using the
  kernel $\KK$ and the initial mean $\mm$. Then, there are constants
  $h_0,C>0$ such that for all  $X$ with $h_{X,\Omega}\le h_0$, all $1\le q\le\infty$ 
    and all $0\le s\le \ell$ with $\ell$ given in Theorem
    \ref{thm:sampling} the following bounds hold.
  \begin{enumerate}
  \item
    If $\beta\ge \tau$ then 
    \begin{eqnarray*}
      \lefteqn{\hspace*{-5ex}\E\left[\|\vv-\mm_{N,\lambda}\|_{W_q^s(\cD)^d}\right] \le C\left( h_{X,\cD}^{\tau-s -
          d(1/2-1/q)_+}
        +\tlg{\sqrt{\lambda}} h_{X,\cD}^{d/\gamma-s}\right)\is{\cN}_\tau(\vv,\mm)}\\
      &&\mbox{}+  C\left(
      \tlg{\frac{1}{\sqrt{\lambda}}} h_{X,\cD}^{\tau-s-d(1/2-1/q)_+ }  +
      h_{X,\cD}^{d/\gamma-s}\right) \is{\E\left[\left(\sum_{j=1}^N \|\eepsilon_j\|_2^2\right)^{1/2}\right]}.
    \end{eqnarray*}
  \item
    If $\beta<\tau$ then
    \begin{eqnarray*}
      \lefteqn{ \E\left[\|\vv-\mm_{N,\lambda}\|_{W_q^s(\cD)^d}\right] \le C h_{X,\cD}^{\beta-s -
          d(1/2-1/q)_+} \rho_{X,\cD}^{\tau-\beta} \is{\cN}_\beta(\vv,\mm)}\\
      &&\mbox{}
      + C \tlg{\sqrt{\lambda}}  h_{X,\cD}^{\beta-\tau +
        d/\gamma-s}\rho_{X,\cD}^{\tau-\beta}  \is{\cN}_\beta(\vv,\mm)\\
      &&\mbox{}+  C\left( \tlg{\frac{1}{\sqrt{\lambda}}}
      h_{X,\cD}^{\tau-s-d(1/2-1/q)_+ } +
      h_{X,\cD}^{d/\gamma-s}\right) \is{\E\left[\left(\sum_{j=1}^N \|\eepsilon_j\|_2^2\right)^{1/2}\right]},
    \end{eqnarray*}
where \is{the} expectation is taking with
respect to the distribution of the errors $\eepsilon_j$, and $\gamma=\max\{2,q\}$.
  \end{enumerate}
\end{theorem}
\begin{proof}
In this situation, we have $\mm_{N,\lambda}=\mm-\is{\QQ_{X,\lambda}}\mm +
\is{\QQ_{X,\lambda}}\vv + \is{\QQ_{X, \lambda}}\cEE$, where again
$\cEE=(\eepsilon_1,\ldots,\ee_N)$. Thus, in addition to the 
error \is{bound} in Theorem \ref{thm:approximationnonoise} we have the
additional term $\|\is{\QQ_{X,\lambda}}\cEE\|_{\is{W_q^s(\cD)^d}}$ that we need
to bound. \is{On the one hand Lemma \ref{lem:standardbounds} for the case
$\vv=\00$, together with Theorem
\ref{thm:nativespace}, yields} the bound
$\|\is{\QQ_{X,\lambda}}\cEE\|_{H^\tau(\cD)^d} \le \frac{C}{\sqrt{\lambda}}\|\cEE\|_{\tlg{\ell_2}}$,
and, on the other hand, \is{also} the bound  
\[
\|\tlg{\QQ_{X,\lambda}\cEE|X}\|_{\tlg{\ell_2}}\le \|\cEE -
\is{\QQ_{X,\lambda}}\cEE|X\|_{\tlg{\ell_2}} +\|\cEE\|_{\tlg{\ell_2}} \le 2\|\cEE\|_{\tlg{\ell_2}}.
\]
Thus the sampling inequality yields in this case
\begin{eqnarray*}
  \|\is{\QQ_{X,\lambda}} \cEE\|_\is{{W_q^s(\cD)^d}} &\le& C
  h^{\tau-s-d(1/2-1/q)_+ }\|\is{\QQ_{X,\lambda}} \cEE\|_{H^\tau(\cD)^d}\\
  & & \mbox{} + C
  h^{d/\gamma-s} \|\is{\QQ_{X,\lambda}}\cEE\|_{\ell_2(X)^d}\\
  &\le& C
  h^{\tau-s-d(1/2-1/q)_+ } \frac{1}{\sqrt{\lambda}}\|\cEE\|_{\tlg{\ell_2}}\\
  & & \mbox{} + C
  h^{d/\gamma-s} \|\cEE\|_{\tlg{\ell_2}}.
\end{eqnarray*}
Together with the error bound from Theorem
\ref{thm:approximationnonoise}, this yields the stated error estimate.
\end{proof}

\begin{remark}
  The above theorem is the  result for vector-valued divergence- or
  rotation-free functions corresponding to 
  \cite[Theorem 4]{Wynne-etal-21-1} for the scalar-valued
  case. Though, we  essentially 
  achieve the same approximation orders, the proof of
  \cite[Theorem 4]{Wynne-etal-21-1} seems problematic, as in
  formula (18) the authors essentially use a bound of the form
  \[
  \|I_X f\|_{H^\tau(\cD)} \le C q_X^{\beta-\tau} \|f\|_{H^\beta(\cD)},
  \]
  which they claim can be derived from the proof of \cite[Theorem
    4.2]{Narcowich-etal-06-1}. We do not believe that this is the
  case, as the above bound \is{would} imply the {\em inverse estimate}
  \[
  \|s\|_{H^\tau(\cD)} \le C q_X^{\beta-\tau} \|s\|_{H^\beta(\cD)},
  \qquad s\in V_X,
  \]
  which is still an open question for general bounded domains. Nonetheless,
  the results of \cite[Theorem
    4]{Wynne-etal-21-1} are correct, as our proof can essentially be
  given also in the scalar-valued case.
\end{remark}

\begin{remark}
Note that in the error bound  of the latter theorem, there are
apparently competing terms on how to choose the smoothing parameter $\lambda$. To
cast more light on this, we look at the specific situation of $q=2$
and $s=0$. If we also assume quasi-uniform data sets then
the error bound for the case $\beta<\tau$,
becomes
\begin{eqnarray*}
  \E\left[\|\vv-\mm_{N,\lambda}\|_{L_2(\is{\cD)^d}}\right] & \le &
  C \left( h_{X,\cD}^\beta + \tlg{\sqrt{\lambda}} h_{X,\cD}^{\beta-\tau+d/2} \right)
  \is{\cN}_{\beta}(\vv,\mm)\\  & & \mbox{} + C \left(\frac{1}{\sqrt{\lambda}} h_{X,\cD}^\tau  +
  h_{X,\cD}^{d/2}\right)\is{\E\left[\left(\sum_{j=1}^N \|\eepsilon_j\|_2^2\right)^{1/2}\right]}.
  \end{eqnarray*}
Thus, the competing terms can for example be balanced by choosing
the smoothing parameter as $\tlg{\sqrt{\lambda}}=h_{X,\cD}^{\tau-\beta/2}$, 
which yields
\begin{eqnarray*}
\E\left[\|\vv-\mm_{N,\lambda}\|_{L_2\is{(\cD)^d}}\right]
&\le& C  h_{X,\cD}^{\min\{\beta, (\beta+d)/2\}} \is{\cN}_{\beta}(\vv,\mm)\\
& & \mbox{}
+ h_{X,\cD}^{\min\{\beta/2, d/2\}} \is{\E\left[\left(\sum_{j=1}^N \|\eepsilon_j\|_2^2\right)^{1/2}\right]}.
\end{eqnarray*}
\end{remark}

For our final result, we assume that each noise $\eepsilon_j$ has \is{an independent} normal distribution
with variance $\sigma^2 I_{d}$ and that we use the
approximation operator with the correct parameter $\lambda=\sigma$.

\begin{corollary} Let the assumptions of Theorem
  \ref{thm:approximationnoise} hold. Assume that $\eepsilon_j\sim
  \cN(\00, \sigma^2 I_{d})$, \is{$1\le j\le N$}.
\begin{enumerate}
  \item
    If $\beta\ge \tau$ then 
    \begin{eqnarray*}
      \lefteqn{\hspace*{-5ex}\E\left[\|\vv-\mm_{N,\sigma}\|_{W_q^s(\cD)^d}\right] \le C\left( h_{X,\cD}^{\tau-s -
          d(1/2-1/q)_+}
        +\sigma h_{X,\cD}^{d/\gamma-s}\right)\is{\cN}_\tau(\vv,\mm)}\\
      &&\mbox{}+  C\rho_{X,\cD}^{d/2}\left(
      h_{X,\cD}^{\tau-s-d(1/2-1/q)_+ -d/2}  +
      \sigma h_{X,\cD}^{d/\gamma-s-d/2}\right).
    \end{eqnarray*}
  \item
    If $\beta<\tau$ then
    \begin{eqnarray*}
      \lefteqn{ \E\left[\|\vv-\mm_{N,\sigma}\|_{W_q^s(\cD)^d}\right] \le C h_{X,\cD}^{\beta-s -
          d(1/2-1/q)_+} \rho_{X,\cD}^{\tau-\beta} \is{\cN}_\beta(\vv,\mm)}\\
      &&\mbox{}
      + C \sigma h_{X,\cD}^{\beta-\tau +
        d/\gamma-s}\rho_{X,\cD}^{\tau-\beta}  \is{\cN}_\beta(\vv,\mm)\\
    &&\mbox{}+  C\rho_{X,\cD}^{d/2}\left(
      h_{X,\cD}^{\tau-s-d(1/2-1/q)_+ -d/2}  +\sigma
      h_{X,\cD}^{d/\gamma-s-d/2}\right).
    \end{eqnarray*}
    \end{enumerate}
\end{corollary}
\begin{proof}
This immediately follows using the estimates from Theorem
\ref{thm:approximationnoise},  together with $N\le C q_X^{-d}=C h_{X,\cD}^{-d}
\rho_{X,\cD}^d$ and
\[
\is{\E\left[\left(\sum_{j=1}^N \|\eepsilon_j\|_2^2\right)^{1/2}\right]} \le\is{ \left(\E\left[
\sum_{j=1}^N \|\eepsilon_j\|_2^2\right]\right)^{1/2}} = \sqrt{N} \sigma \; \is{\le \; h_{X,\cD}^{-d/2} \rho_{X,\cD}^{d/2} \sigma}.
\]
\end{proof}

\section*{Acknowledgments}

{Quoc Thong Le Gia and Ian Sloan are  grateful to the
  Australian Research Council for support under grant DP180100506. Holger Wendland is grateful for the support and hospitality of the Sydney Mathematical
  Research Institute (SMRI) during his visit in August/September
  2023.}

\end{document}